\documentclass[a4paper]{article}
\usepackage{amsmath,amsfonts,amsthm,amssymb}
\usepackage[shortlabels]{enumitem}
\usepackage{natbib}
\usepackage[dvipsnames]{xcolor}
\usepackage{hyperref}
\hypersetup{
    colorlinks,citecolor=blue,urlcolor=blue,linkcolor=blue,linktocpage=true
}

\usepackage{xr}
\newtheorem{theorem}{Theorem}
\newtheorem{remark}{Remark}

\title{A central limit theorem for a sequence of conditionally
  centered random fields}

\author{Abdollah Jalilian\\ Department of Statistics, Razi
University, Iran\\ Lancaster Ecology and Epidemiology Group,\\
Lancaster University, United Kingdom\\\texttt{stat4aj@gmail.com}
\and Arnaud Poinas\\ Laboratoire de Mathematiques et Applications,\\
University of Poitiers,
France\\\texttt{arnaud.poinas@univ-poitiers.fr}
\and Ganggang Xu\\
Department of Management Science,\\ University of Miami, United States
of America\\\texttt{gangxu@bus.miami.edu}
\and Rasmus Waagepetersen\\ Department of Mathematical
Sciences, \\ Aalborg University, Denmark\\\texttt{rw@math.aau.dk}}%

\def\H{\mathcal H}
\def\R{\mathbb R}
\def\Z{\mathbb{Z}}
\def\EE{\mathbb E}
\def\Var{{\mathbb{V}\mathrm{ar}}}
\def\Cov{{\mathbb{C}\mathrm{ov}}}
\def\Corr{{\mathbb{C}\mathrm{orr}}}

\def\T{\mathsf{T}}
\def\N{\mathbb N}

\def\1{\mathbf{1}}
\def\X{\mathbf{X}}

\newcommand{\Ekn}[1]{ E_{{#1}} }

\newcommand{\mkn}{m_n}
\newcommand{\lkn}{l_K}

\def\dd{\mathrm{d}}

\def\i{\mathbf{i}}
\def\j{\mathbf{j}}
\def\l{\mathbf{l}}
\def\k{\mathbf{k}}

\def\D{\mathcal{D}}

\reversemarginpar

\newcommand{\e}{{\text e}}

\newcommand{\bsl}{\baselineskip}
\newcommand{\rwrev}[1]{{#1}}

\newcommand{\apMOD}[1]{{#1}}

\begin{document}
\maketitle

\begin{abstract}
A central limit theorem is established for a sum of random variables
belonging to a sequence of random fields. The fields are assumed to
have zero mean conditional on the past history and to satisfy certain
conditional $\alpha$-mixing conditions in space or time. Exploiting
conditional centering and the space-time structure, the limiting
normal distribution is obtained for increasing spatial domain,
increasing length of the sequence, or both of these. The theorem is very well
  suited for establishing asymptotic normality in the context of unbiased
  estimating function inference for a wide range of space-time
  processes. This is pertinent given the abundance of space-time
    data. Two examples demonstrate the applicability of the theorem.

\textbf{MSC2020 subject classifications:} 60F05; 60G55; 60G60; 62M30
\end{abstract}

\section{Introduction}

While the literature on central limit theorems for
spatial stochastic processes is extensive \cite[][to mention a few]{bolthausen:82,jensen1994,guyon:95,comets:janzura:98,karacsony:06,biscio:waagepetersen:18}, the literature on central limit theorems
for spatio-temporal processes seems much less developed. The theoretical foundation for statistical analysis of
  space-time data is thus incomplete which is unfortunate since such
  data are increasingly common.
The objective of this paper is to derive a
central limit theorem for a sequence of random fields that satisfy a
certain conditional centering condition. This condition is
motivated by applications of estimating functions in
space-time statistics (see Section~\ref{sec:applications})
where the estimating function based on a random field at
one time point is unbiased (has zero mean expectation at the true
parameter value) conditional on the past.

Our central limit theorem is based on increasing spatial
domain, increasing time horizon or both but does not require stationarity in space nor in time. To prove the theorem we use a unified framework that more specifically covers both a fixed length
sequence of random fields each with increasing spatial
domain, a sequence of random fields with unconstrained
time horizon and fixed spatial domain, and a combination of these
asymptotic settings.
The proof requires to bound
a large variety of covariances. In the case of spatial random fields,
this is typically done by imposing appropriate weak spatial dependence
conditions like $m$-dependence \cite[e.g.][]{leonenko:75,Heinrich2013} or $\alpha$-mixing \citep{bolthausen:82}. \cite{leonenko:75} considers
a central limit theorem for a sequence of statistics each obtained as
a sum of variables for a random field in a sequence of increasing domain $m$-dependent random
fields. Instead of using mixing conditions as an intermediate step to
bound covariances, \cite{doukhan:lang:16} propose a notion of weak
dependence where they impose decay
conditions directly on covariances. \cite{jensen1994} and \cite{comets:janzura:98} use a different \rwrev{strategy} for spatial random fields of conditionally centered random variables that have expectation zero given all other
variables. By exploiting the conditional centering, they are able to
dispense with spatial mixing assumptions.

Different from \cite{jensen1994} and \cite{comets:janzura:98}, we assume
conditional centering given the past history in time, which is useful for eliminating
covariances involving variables at distinct time points. This is related
to the key property of uncorrelated increments for martingales. In
addition we exploit the space-time structure and introduce a novel
assumption of $\alpha$-mixing for a random field at a given point in time
conditional on the past history. This is natural to consider when the distribution of a space-time process is specified in
terms of a sequence of conditional distributions given the past, see
Section~\ref{sec:applications}. Referring again to martingale terminology, this
type of $\alpha$-mixing is convenient since it essentially only
pertains the spatial mixing properties of innovations.
  We note that our notion of conditional
$\alpha$-mixing is very different from the notion of conditional
$m$-dependence in \cite{biscio:svane:22} for spatial point processes
where conditioning is on an underlying random field.

To give a concise statement of the various conditions needed, we
introduce the novel concept of a space-time filtration which is an
increasing sequence of $\sigma$-algebras that are indexed both by time
points and
subsets of space.

\section{Examples of applications}\label{sec:applications}

To motivate our central limit we provide some illustrative examples of applications where we try to strike a balance between
accessibility and practical relevance. The targets for our central
limit are so-called score functions. These play a key role in
statistical inference where parameter estimates are obtained by
solving estimating equations where score functions involving unknown
parameters are equated to zero. \rwrev{We use the first example as a simple illustration of the kind of localization and mixing conditions that we impose for our central limit theorem. }

\subsection{\rwrev{A vector autoregression of order one}}\label{sec:var}

\rwrev{
Let $X_{0}=(X_{0}(\l))_{\l \in \Z}$ denote an initial
condition and
let  $\nu_{k}=(\nu_k(\l))_{\l \in \Z}$, $k \ge 1$, denote a sequence
of independent and identically distributed zero-mean stationary processes. Then for $n=1,2,\ldots$ we for consistency of notation
let $X_{0,n}(\l)=X_0(\l)$, $\l \in \Z$, and define for $k>0$,
\[ X_{k,n}(\l)= a \bar X_{k-1,n}(\l)+\nu_k(\l), \quad \l=1,\ldots,n, \]
where $\bar X_{k-1,n}(\l)$ is an average of $X_{k-1,n}(\l)$ and its immediate neighbours and $a \in \R$ is
an autogression coefficient. More specifically, $\bar X_{k-1,n}(1)=[X_{k-1,n}(1)+X_{k-1,n}(2)]/2$, $\bar X_{k-1,n}(\l)=[X_{k-1,n}(\l-1)+X_{k-1,n}(\l)+X_{k-1,n}(\l+1)]/3$ for $1<\l<n$, and $\bar X_{k-1,n}(n)=[X_{k-1,n}(n-1)+X_{k-1,n}(n)]/2$. 
Suppose the $X_{k,n}(\l)$ are observed for $k=0,\ldots,K$ and $\l \in \D_n=\{1,\ldots,n\}$. We may then estimate $a$ by solving $T_{K,n}=0$ with respect to $a$ where $T_{K,n}$ is the least squares score function
\[ T_{K,n}=\sum_{k=1}^K \sum_{\l=1}^n \bar X_{k-1,n}(\l) [X_{k,n}(\l)-a \bar
  X_{k-1,n}(\l)].\]

Despite the simple formulation of this model, the space-time
correlation structure is rather complicated. For a given
time point $k$, spatial dependence between two variables $X_{k,n}(\l)$
and $X_{k,n}(\j)$ arise to begin with from the spatial dependence in $X_0$ and the
$\nu_l$ processes for any  $1 \le l \le k$. In addition, spatial
dependence may be due to common previous terms $X_0(\k)$, $\nu_l(\k)$, $1 \le l <k$, $\k \in
\Z$ appearing in the sums defining $X_{k,n}(\l)$ and
$X_{k,n}(\j)$.

The model is closely related to the network vector
autoregression (NAR) introduced in \cite{zhu2017network} for large-scale (large $n$) networks, and its numerous variants, e.g.,  \cite{zhu2023simultaneous} and \cite{liu2024two}. It is
simpler than the NAR in the sense that we consider a very specific neighbour
structure and do not consider covariates. We on the
other hand  allow the innovation terms
$\nu_k(\l)$, $\l \in \Z$, to be correlated while \cite{zhu2017network} consider independent innovations. Using our central limit
theorem it turns out to be relatively straightforward to obtain
asymptotic normality for $T_{K,n}$ and our theorem seems useful as
well for more general NAR models with correlated innovations.}

\subsection{A discrete time spatial birth-death point process}\label{sec:iteratedpointprocess}

Consider a sequence of point processes
$X_0,\ldots,X_K$ and spatial covariates $Z_0,\ldots,Z_K$ on $\R^d$,
where the latter are considered non-random. For each $k=0,\ldots,K$, we observe $X_k \cap
W$ and $Z_k(v)$, $v \in W,$ where $W \subset \R^d$ is a bounded observation window. We consider a model for the $X_k$
corresponding to a discrete-time spatial birth-death process that mimics the seasonal behaviour of a plant community.

Given $X_{k-1}$ and $Z_{k-1}$, $X_k = S_{k} \cup B_{k} \cup I_k$ where $B_{k} = \cup_{u \in X_{k-1}} Y^k_u$ is a
union of offspring clusters $Y^k_u$ given by independent Poisson
processes. For each $u \in X_{k-1}$,  $Y_u^k$ has an intensity function given by $\alpha k_u(\cdot)$ where $\alpha >0$ and $k_u(\cdot)$ is the uniform density on the disc  $b(u,\omega)$ with center $u$ and
radius $\omega$. The intensity function depends on $X_{k-1}$ only through $u$. The point process $S_{k}$ represents plants from $X_{k-1}$ that survived to time $k$ and is
given by $S_{k} = \{ u \in X_{k-1} | R_u =0\}$, where the $R_u$ are
Bernouilli variables conditionally independent given $X_{k-1}$ and
$Z_{k-1}$ and with `survival' probability
$\mathbb{P}(R_u=0|X_{k-1},Z_{k-1})$ depending on $X_{k-1}$ and $Z_{k-1}$
only through $Z_{k-1}(u)$. The `immigrant' point processes
  $I_k$, $k \ge 1$, are independent Poisson processes of intensity
  $\rho \ge 0$. Given $X_{k-1}$, the process $B_k \cup I_k$ of new trees in generation $k$ is a Poisson point process with intensity function given by $\rho+\alpha \sum_{u \in X_{k-1}} k_u(v)$, $v \in W$.

Suppose parametric models are imposed for the conditional distributions of $B_k \cup I_k$ and $S_k$ given $Z_{k-1}$ and $X_{k-1}$. For each $k$ we let $T_{k}^{W}$ denote the concatenation of likelihood score functions for the conditional distributions of $B_k  \cup I_k$ and $S_k$ . Then $T_{k}^{W}$ can be additively
decomposed as
\[ T_{k}^{W} = \sum_{\l \in \D} \Ekn{k}(\l),
\]
where $\D=\{\l \in \Z^d: c(\l) \cap W \neq \emptyset \}$, $c(\l)$
is the unit cube centered at $\l$ and the likelihood score component
$E_{k}(\l)$ depends on
$Z_{k-1},X_{k-1},B_k\cup I_k,S_k$ only through $Z_{k-1}(v)$, $v \in
c(\l) \cap W$ and through the intersections of $X_{k-1},B_k \cup I_k,$
and $S_k$ with $c(\l) \cap W$. Due to boundary effects, the definition of
the variables $\Ekn{k}(\l)$ may depend on $W$ (when $c(\l)$ is not
entirely contained in $W$). Therefore we will later on, when
considering a sequence of expanding windows $W_n$, add the index $n$
and write $\Ekn{k,n}(\l)$.

Parametric inference is based on the accumulated score function
\[ T^{W} = \sum_{k=1}^K T_k^W \]
which has zero mean (is unbiased) since each conditional score
function has zero conditional mean by the first Bartlett identity. For asymptotic inference it is essential to establish asymptotic normality of $T^W$. Here relevant asymptotic regimes can be expanding
time horizon $K \rightarrow \infty$, expanding observation window $W$
or a combination of these.

The case of increasing window asymptotics
is pertinent for rain forest point pattern data sets covering large
study regions, \rwrev{see \cite{jalilian2024compositelikelihoodinferencespacetime}}. Here tree positions are registered in
censuses conducted a moderate number of times. For example, for the
Barro Colorado Island data set \citep{condit:etal:19}, censuses are
available at just eight time points which does not justify large $K$
asymptotics.

\subsection{\rwrev{Graphical} autoregressive model}\label{sec:spacetimecar}

In this section, we consider a spatio-temporal autoregressive model. The model can also be viewed
as a multivariate time series with graphical interactions at each time
instance \citep{dahlhaus2000graphical}.
Let $\D \subset\R^{d}$ be a finite lattice and $G=(\D, \mathcal{E})$ be an undirected graph with vertices in $\D$ and edges $\mathcal{E}$ where
$\l \neq \j$ in $\D$ are connected if $\{\l,\j\} \in
\mathcal{E}$.
 For example, elements of $\D$ may represent sub-regions (provinces, counties or municipalities) of a specific geographical region and edges in $\mathcal{E}$ determine interconnections between these sub-regions due to {e.g.\ adjacency or transportation routes}.

Let $X_{k}=\left( X_{k}(\l) \right)_{\l\in\D}$, $k\in\Z_{+}=\{0,1,2,\ldots\}$, be a sequence of random fields where for any $k\in\Z_{+}$ and $\l\in\D$, we assume that given the history $X_{0},\ldots,X_{k-1}$ and the current neighbourhood $\rwrev{X_k(-\l)=\{X_{k}(\j)\}_{\j \in \D \setminus \{\l\}}}$,  $X_{k}(\l)$ follows the conditional distribution
\begin{equation}\label{eq:condspec}
  X_{k}(\l) \big| X_{0},\ldots,X_{k-1},X_k(-\rwrev{\l}) 
\sim
  \mathcal{N}\left( \beta \xi_{k-1}^{\mathrm{temp}}(\l)+ \gamma \xi_{k}^{\mathrm{spat}}(\l) ,1/a(\l) \right),
\end{equation}
where for some $r \ge 1$,
\[
  \xi_{k-1}^{\mathrm{temp}}(\l)= \sum_{j=1}^{r} \sum_{\j \in \D}
  b_{j}(\l,\j)X_{k-j}(\j),
\]
with $X_{k-j}(\l)=0$ for $k<j$, and
\[
   \xi_{k}^{\mathrm{spat}}(\l) = \sum_{\j \in \D\backslash\{\l\}}   b_0(\l,\j) \left( X_{k}(\j) - \beta \xi_{k-1}^{\mathrm{temp}}(\j) \right).
\]
Let $\xi_{k-1}^{\mathrm{temp}}=(\xi_{k-1}^{\mathrm{temp}}(\l))_{\l\in \D}$, $\xi_{k}^{\mathrm{spat}}=(\xi_{k}^{\mathrm{spat}}(\l))_{\l\in \D}$, $A = \mathrm{diag}\left(a(\l): \l\in\D \right)$, and $B_{j}=\left[b_{j}(\l,\j) \right]_{\l,\j\in \D}$, $j=0,\ldots,r$, with diagonal entries $b_{0}(\l,\l)=0$ for $B_{0}$.
Then,
\[
  \xi_{k-1}^{\mathrm{temp}} = \sum_{j=1}^{r} B_{j} X_{k-j} \quad \text{ and } \quad \xi_{k}^{\mathrm{spat}} = B_{0} \left(X_{k} - \beta \xi_{k-1}^{\mathrm{temp}} \right).
\]

We assume that $b_0(\l,\j),\ldots, b_{p}(\l, \j)$ and $a_{n}(\l)$ are
known and  $\beta,\gamma\in\R$ are unknown parameters. Based on
  observations of $X_{k}$ for $k=0, \ldots, K$,  a score function for estimation of $\beta$ and $\gamma$ is given by
\[
  T = \sum_{k=1}^{K} \sum_{\l\in \D} E_{k}(\l),
\]
\rwrev{\cite[see Section~5.2 
 of the supplementary material][for
 details]{jalilian:poinas:xu:waagepetersen:23:suppl}} where
\[
  E_{k}(\l) = \varepsilon_{k}(\l) \left(\xi_{k-1}^{\mathrm{temp}}(\l), \xi_{k}^{\mathrm{spat}}(\l)\right)^{\T}
\]
and $\varepsilon_{k}=A(X_{k} -  \beta  \xi_{k-1}^{\mathrm{temp}}-\gamma \xi_{k}^{\mathrm{spat}})$.

We can asympotically approximate the complicated distribution of $T$ by a normal distribution as $K$ tends to infinity.

\section{The space-time central limit theorem}\label{s:clt}

Consider for each $k,n\in\N$ a random field $\Ekn{k,n}=(
\Ekn{k,n}(\l) )_{\l \in \D_n}$, on $\D_n  \subseteq \Z^d$, $d\geq 1$,  where
$\Ekn{k,n}(\l)=\left(E^{(1)}_{k,n}(\l),
  \ldots,E^{(q)}_{k,n}(\l)\right)^{\T},$ $q \ge 1$, and $(\D_n)_{n\in\N}$ is a
sequence of finite observation lattices while $k$ can be viewed as a time
index. \rwrev{As explained in Sections~\ref{sec:var} and \ref{sec:iteratedpointprocess},  the index $n$ for the variables $\Ekn{k,n}(\l)$ may be needed due to
  edge effects when the variables are constructed in terms of an autoregression or from an underlying
  process on a continuous observation window $W_n$.}

Depending on the context (see examples in Section~\ref{sec:applications}), it may or may not be reasonable to
assume properties like stationarity or mixing in time for the sequence
	$\Ekn{k,n}$, $k \in \N$. If the sequence
$\Ekn{k,n}$, $k \in \N$, can not be viewed as a part of a stationary
sequence, it is convenient to consider it conditional on an
initial condition represented by a $\sigma$-algebra $\H_0$, see
Section~\ref{s:notation}.
Defining for $K,n \in \N$,
\begin{equation}\label{eq:TWsum}
T_{K,n} = \sum_{k=1}^K \sum_{\l \in \D_n} \Ekn{k,n}(\l),
\end{equation}
we formulate in Section~\ref{s:cltstate} a central limit theorem for
$\Sigma_{K,n}^{-1/2} T_{K,n}$ with $\Sigma_{K,n} =
\Var_{\H_0}[T_{K,n}]$, that is, the variance conditional on the initial condition $\H_0$. Conditions for the theorem are specified in Section~\ref{s:cond}.

\subsection{\rwrev{Initial remarks in relation to vector autoregression}}\label{sec:simpleintrocond}

\rwrev{
Before going deeply into technical definitions and conditions we
return for a moment to the vector autoregression in Section~\ref{sec:var} to
give an informal motivation and explanation of  various technical
concepts introduced in the next sections.

Letting $E_{k,n}(\l)=\bar X_{k-1,n}(\l) [X_{k,n}(\l)-a \bar X_{k-1,n}(\l)]$
it is clear that the $E_{k,n}(\l)$ are conditionally
centered, $\EE[ E_{k,n}(\l)|
X_{k-1.n}]=0$. For $Y^{k,K,\j}=f[E_{k,n}(\j),E_{k+1}(\j),\ldots,E_{K,n}(\j)]$ given as a function of $E_{k,n}(\j)$ and future variables $E_{l,n}(\j)$, $ k<l \le K$, we consider further
the conditional expectation of $Y^{k,K,\j}$ given the past
$X_0,\nu_1\ldots,\nu_{k-1}$. We will formulate such a conditional expectation in terms of a space-time filtration introduced in the next section. For the vector autoregression  we for example
 define the past $\H_{k-1}$ as the $\sigma$-algebra generated
 by $\X_0$ and $\nu_1,\ldots,\nu_{k-1}$ and
 $\H_{0,k-1,\{\j-K,\ldots,\j+K\}}$ as the smaller $\sigma$-algebra
 generated by $X_0(\l)$ and $\nu_{l}(\l)$ for $0 < l \le k-1$ and
 $\l \in \{\j-K,\ldots,\j+K\}$. Then the conditional expectation is  `localized'
 (or Markovian) in the sense that
 $\EE[Y^{k,K,\j}|\H_{k-1}]=\EE[Y^{k,K,\j}|\H_{0,k-1,\{\j-K,\ldots,\j+K\}}]$.
  Spatial localization is
  important when we need to control spatial dependence between for
  example $\EE[Y^{k,K,\l}|\H_{k-1}]$ and $\EE[Y^{k,K,\j}|\H_{k-1}]$,
  $\l \ne \j$,  using spatial $\alpha$-mixing coefficients (next section). This is because the mixing coefficients are
  concerned with dependence between parts of the stochastic process
  where the parts are defined in terms of subsets of space or time or both.

For
  the toy example, if the correlations $\Corr[\nu_k(\l),\nu_{k}(\l+h)]$
  are sufficiently fast decaying as a function of $h$, then $\nu_k$ is
  $\alpha$-mixing \citep{doukhan:94}. This immediately implies that stochastic processes involving $\nu_k$ and the past are $\alpha$-mixing given the past. For example, $E_{k,n}$ with components $\bar X_{k-1,n}(\l) \nu_k(\l)$ is
  $\alpha$-mixing conditional on $\H_{k-1}$. When we in the next section introduce a mixing coefficient for $\H_{0,k,A}$ and $\H_{0,k,B}$ given $\H_{k-1}$ this essentially involves mixing between $\nu_k(\l), \l \in A$, and $\nu_k(\j)$, $\j \in B$.

For fixed $n$ and increasing $K$ we might instead consider
localization in time. In this case we may  define the past
$\H_{0,k-1,\D_n}$ as the $\sigma$-algebra generated by $X_{n,l}(\j)$,
$0 \le l \le k-1$, $\j \in \D_n$, and  $\H_{k-1,k-1,\D_n}$ as the
$\sigma$-algebra generated by $X_{k-1,n}(\j)$, $\j \in \D_n$, in which case the conditional expectation $\EE[Y^{k,K,\j}|\H_{0,k-1,\D_n}]=\EE[Y^{k,K,\j}|\H_{k-1,k-1,\D_n}]$ is localized in time.

}

\subsection{Notation, space-time filtration and
  $\alpha$-mixing}\label{s:notation}

We define $
\mathrm{d}(x,y) = \max \lbrace |x_i-y_i |: 1\leq i \leq d \rbrace,\  x,y \in \R^d,
$
and, reusing notation,
$
 \mathrm{d}(A,B) = \inf \lbrace \mathrm{d}(x,y): x\in A, \ y\in B
 \rbrace,\  A,B \subseteq \R^d.
$ 
 For a subset $A \subseteq \R^d$ we denote by $\left| A \right|$ the
  cardinality or Lebesgue
measure of $A$. The meaning of $|\cdot|$ and $\mathrm{d}(\cdot,\cdot)$ will
be clear from the context. For $a \in \R$ we use the brief notation
$\EE[\cdots]^a$ for the less ambiguous $\left (\EE[\cdots] \right)^a$.

Conditioning on the past history plays a crucial role for our result. We
represent the history by what we call a space-time filtration
whose definition not only involves time points but also spatial
regions. Let $(\Omega, \mathcal{F}, \mathbb{P})$ be the underlying
probability space on which all random objects in this paper are
defined. Recalling $\Z_{+}=\{0, 1, 2,\ldots\}$, we coin a set $\H=\{
\H_{l,k,A} \mid l \le k \in \Z_{+}, A \subseteq \R^d\}$ of
$\sigma$-algebras $\H_{l,k,A}\subseteq \mathcal{F}$ a space time
filtration if $\H_{l,k,A} \subseteq \H_{l',k',A}$ for $l'\le l \le k
\le k'$ and $\H_{l,k,A} \subseteq \H_{l,k,B}$ for $A \subseteq
B$.
For brevity we let $\H_k=\H_{0,k,\R^d}$. \rwrev{As discussed in Section~\ref{sec:simpleintrocond}
one can think of the $\sigma$-algebras $\H_{l,k,A}$
as being generated by a stochastic process restricted to the space-time domain
$\{l,l+1,\ldots,k\}\times A$. More generally, this could also include
space-time exogeneous variables serving as covariates in a
space-time model. This in turn makes the $\Ekn{k,n}(\l)$
measurable when defined in terms of the underlying processes. We consider further examples in
Sections~\ref{sec:verifyiteratedpointprocess} and
\ref{sec:verifyspacetimecar}}.

For any $\sigma$-algebras $\mathcal{G}_0, \mathcal{G}_1,
\mathcal{G}_2\subseteq\mathcal{F}$, the conditional $\alpha$-mixing
coefficient of $\mathcal{G}_1$ and $\mathcal{G}_2$ given
$\mathcal{G}_0$ is defined in \cite{rao:09}  as
\[
  \alpha_{\mathcal{G}_0}(\mathcal{G}_1, \mathcal{G}_2) = \sup_{\substack{F_1\in\mathcal{G}_1\\ F_2\in\mathcal{G}_2}} \big|\mathbb{P}_{\mathcal{G}_0}(F_1\cap F_2) - \mathbb{P}_{\mathcal{G}_0}(F_1) \mathbb{P}_{\mathcal{G}_0}(F_2)\big|,
\]
where $\mathbb{P}_{\mathcal{G}_0}(F) = \EE_{\mathcal{G}_0}[I_{F}]$,
$\EE_{\mathcal{G}_0}$ denotes the conditional expectation with respect
to $\mathcal{G}_0$, $I_{F}(\omega)=I[\omega \in F]$ and $I[\,\cdot\,]$
is the indicator function. Similarly, $\Var_{\mathcal{G}_0}$
  and $\Cov_{\mathcal{G}_0}$ denote variance and covariance
  conditional on $\mathcal{G}_0$. Note that $0\leq
  \alpha_{\mathcal{G}_0}(\mathcal{G}_1, \mathcal{G}_2) \leq 1/4$ is a
  $\mathcal{G}_0$-measurable random variable.

For $m,c_1,c_2\geq 0$ and $K\in\N$ we
define a spatial mixing coefficient
\begin{align}\label{eq:alpha}
  \alpha_{K,c_1,c_2}(m) &= \sup_{1 \le k \le K} \sup_{\substack{A,B\subseteq \R^d \\ |A| \leq c_1, |B|\leq c_2 \\  \dd(A, B) \geq m}} \alpha_{\H_{k-1}}(\H_{0,k,A}, \H_{0,k, B}).
\end{align}
Moreover, for $m, r \ge 0$, $n \in \N$, and $D \subseteq \R^d$,
we define 
 the time mixing coefficient
\begin{align}\label{eq:alphatime}
 \rwrev{ \alpha_{r}(m)} &\rwrev{= \sup_{\substack{l, k \in \Z_{+} \\ l-k \ge m }}
    \alpha_{\H_{k-1}}(\H_{k}, \H_{(l-r+1)^{+},l,\R^d})}
\end{align}
where $(x)^+=\max\{x,0\}$. Finally, we combine \eqref{eq:alpha} and \eqref{eq:alphatime} to
obtain for  $m, m', c_1,c_2, r\geq 0$ \apMOD{and $K \in \N$} the space-time mixing coefficient
\apMOD{\begin{align}\label{eq:alphaspacetime}
  \alpha_{r; K,c_1,c_2}(m; m') &= \sup_{\substack{l, k \leq K \\ l-k \ge m }}\sup_{\substack{A,B\subseteq \R^d\\ |A| \leq c_1, |B|\leq c_2 \\  \dd(A, B) \geq m'}} \alpha_{\H_{k-1}}(\H_{0,k,A}, \H_{\rwrev{(l-r+1)^{+}},l,B}).
\end{align}}
This extends the previous mixing coefficients since \eqref{eq:alpha}
is \apMOD{$\alpha_{\infty; K,c_1,c_2}(\rwrev{0;m})$} \rwrev{if we restrict $l\ge k$ to $l=k$} and \eqref{eq:alphatime}
corresponds to $\alpha_{\rwrev{r}; \infty,\infty,\infty}(\rwrev{m;0})$.

\begin{remark}
The definitions \eqref{eq:alpha} and \eqref{eq:alphatime} involve in a standard way
$\sup$ over
space lags as well as over sizes of spatial regions or over time
lags. The remarkable aspect of \eqref{eq:alpha} is that it just
involves conditional $\alpha$-mixing for events up to time $k$
conditional on the history up to the immediate past at time
$k-1$. This is extremely convenient when considering the common case
of processes  specified sequentially by the distribution of the
  process at time $k$, $k \ge 1$, given the history up to time $k-1$,
  see also Section~\ref{sec:alphamixing}. For \eqref{eq:alphatime} we
  similarly only need to consider the development of the of process
  from time $k$ up to $l$ given the history up to time
  $k-1$. The coefficient \eqref{eq:alphaspacetime} \rwrev{follows by combining} \eqref{eq:alpha} and \eqref{eq:alphatime}.
\end{remark}

\subsection{Conditions}\label{s:cond}

This section specifies and discusses conditions needed for the central limit theorem. We initially assume the existence of a
space-time filtration $\H$ so that the following conditional centering
condition holds:
\begin{enumerate}[label=($\mathcal{A}$\arabic*)]
\setcounter{enumi}{0}
\item \label{cond:condcentered} for any $k,n\in \N$, the random vectors $\Ekn{k,n}(\l)$, $\l\in\D_n$, are conditionally centered:
\begin{equation*}
  \EE_{\H_{k-1}} [ \Ekn{k,n}(\l) ] =0.
\end{equation*}
\end{enumerate}
We further need certain space or time `localization' \rwrev{(or
  Markovian)} properties:
\begin{enumerate}[label=($\mathcal{A}$\arabic*)]
\setcounter{enumi}{1}
\item \label{cond:localizationgeneral}
Consider $n,k \in \N$, $\l\in\Z^{d}$,
$\mathcal{I} \subset \Z^d$, $\mathcal{T}\subseteq\{k, k+1, \ldots, K\}$.
  Let $f:\R^{q|\mathcal{I}||\mathcal{T}|} \to \R$ be any measurable
  function and define $Y_{n}^{(\mathcal{T}, \mathcal{I})} = f \left(
    (E_{k',n}(\l))_{(k', \l) \in \mathcal{T}\times \mathcal{I}}
  \right)$.
\begin{enumerate}[label=($\mathcal{A}$\arabic{enumi}-\arabic*)]
\item \label{cond:spacelocalized} there exists a distance $R>0$ so that $\Ekn{k,n}(\l)$ and
  $\EE_{\rwrev{\H_{k-1}}}[ Y_{n}^{(\mathcal{T}, \mathcal{I})} ]$ are measurable
  with respect to respectively $\H_{0,k,C(\l)}$ and
  $\H_{0, k-1,\cup_{\l \in \mathcal{I}} C(\l)}$, where $C(\l)$ denotes a cube of
  sidelength  $R$ centered at $\l$.
\item \label{cond:timelocalized} there exists a lag
  $r \in \N$ so that $\Ekn{k,n}(\l)$ and  $\EE_{\H_{k-1}}[ Y_{n}^{(\mathcal{T}, \mathcal{I})} ]$ are
  measurable with respect to  respectively \rwrev{$\H_{k}$  and
  $\H_{(k-r)^+,k-1,\R^d}$}.
\item \label{cond:spacetimelocalized} there exists
  a distance $R>0$ and a lag $r \in \N$ so that $\Ekn{k,n}(\l)$ and
  $\EE_{\H_{i-1}}[ Y_{n}^{(\mathcal{T}, \mathcal{I})} ]$, for all $i$ in $\{1,\ldots, k\}$, are measurable
  with respect to respectively $\H_{0,k,C(\l)}$ and \rwrev{$\H_{(i-r)^+, i-1,((k-i)\cup_{\l \in \mathcal{I}} C(\l))}$}, where $C(\l)$ denotes a cube of
  sidelength  $R$ centered at $\l$.
\end{enumerate}
\end{enumerate}
Condition~\ref{cond:condcentered} implies that
$\EE_{\H_{k'}}[\Ekn{k,n}(\l)]=0$ for any $k' <k$. Combined with
condition~\ref{cond:localizationgeneral} it further holds for
any $0 \le k'< k < l$, and $\l,\j\in\Z^d$,
\begin{align*}
  \Cov_{\H_{k'}}[ \Ekn{k,n}(\l), \Ekn{l,n}(\j) ]
  = \EE_{\H_{k'}}\left[ \EE_{\H_{l-1}}[ \Ekn{k,n}(\l) \Ekn{l,n}^{\T}(\j)]  \right] = 0.
\end{align*}
Moreover, it holds that $\EE_{\H_{0}}[T_{K,n}]=0$ and
\[
  \Sigma_{K,n} = 
  \sum_{k=1}^{K} \Var_{\H_{0}}\Big[\sum_{\l \in \D_n} \Ekn{k,n}(\l)  \Big] = \sum_{k=1}^{K}\sum_{\l,\j\in \D_n} \EE_{\H_{0}}\left[ \Cov_{\H_{k-1}}[\Ekn{k,n}(\l), \Ekn{k,n}(\j) ] \right].
\]
Next, we introduce the following conditional $\alpha$-mixing
conditions for the filtration $\H$:
\begin{enumerate}[label=($\mathcal{A}$\arabic*)]
\setcounter{enumi}{2}
\item there exists a real number $p>2$ such
  that \label{cond:conditionalmdependent} one of the following holds
  almost surely (where $r$ is defined in $(\mathcal{A}2)$).
\begin{enumerate}[label=($\mathcal{A}$\arabic{enumi}-\arabic*)]
\item for a fixed $K\in\N$, there exists a function
  $\bar{\alpha}:\R_{+}\to \R_{+}$ such that
\[
  \EE_{\H_{0}}\left[  \alpha_{K,2R^2, \infty}(m) \right] ^{\frac{p-2}{p}} \leq \bar{\alpha}(m)
\]
and $\bar{\alpha}(m)=O(m^{\eta})$, where
$\eta+d<0$.
    \label{cond:spatial conditional independence}
\item for a fixed $n\in\N$, there exists a function $\tilde{\alpha}:
  \R_{+}\to \R_{+}$ such that
\[
  \EE_{\H_{0}}\left[ \alpha_{\rwrev{r}}(m)\right] ^{\frac{p-2}{p}} \leq \tilde{\alpha}(m)
\]
and $\tilde{\alpha}(m)=O(m^{\psi})$, where $\psi< -1$.
      \label{cond:temporal conditional independence}
\item \label{cond:space-time conditional independence}
there exists a function
  $\tilde{\alpha}:\R_{+}\times \R_{+}\rightarrow \R_{+}$ such that
\[
  \sup_{k\in\N}\EE_{\H_{0}}\left[  \apMOD{\alpha_{r;k,2R^2,\infty}}(m; m') \right] ^{\frac{p-2}{p}} \leq \tilde{\alpha}(m; m')
\]
and  $\tilde{\alpha}(m;m')=O(m^{\psi}m'^{\eta})$, where
$\psi< -1$ and $\eta+d<0$. In this case we write
  $\bar\alpha(m)=\tilde\alpha(0; m)$ in analogy to \ref{cond:spatial conditional independence}.
\end{enumerate}
\end{enumerate}
Finally, to obtain a central limit theorem for $\Sigma_{K,n}^{-1/2}T_{K,n}$, we
  need the following sets of conditions:
\begin{enumerate}[label=($\mathcal{A}$\arabic*)]
\setcounter{enumi}{3}
\item there exist an $\epsilon   \geq 3(p-2)$ and a constant
  $0<M_{\epsilon}<\infty$ such that for
\[
  d_{K, n} = \sup_{1\leq  k\leq K} 
  \sup_{\l \in \D_n} \sup_{i=1,\ldots,q} \EE_{{\H_0}}\left[ \left|  \Ekn{k,n}^{(i)}(\l) \right|^{6 + \epsilon} \right],
\]
one of the following holds almost surely
\label{cond:boundedmoments}
\begin{enumerate}[label=($\mathcal{A}$\arabic{enumi}-\arabic*)]
\item  $\sup_{n \in \N} d_{K,n} \leq M_{\epsilon}^{6+\epsilon}$.
\label{cond:boundedmoment1}
\item $\sup_{K \in \N} d_{K,n} \leq M_{\epsilon}^{6+\epsilon}$.
\label{cond:boundedmoment2}
\item $\sup_{K\in\N}\sup_{n \in \N} d_{K,n} \leq
   M_{\epsilon}^{6+\epsilon}$.
 \label{cond:boundedmoment3}
\end{enumerate}
\end{enumerate}
and, with $\lambda_{\min}(M)$ being the smallest eigen value of a
  symmetric matrix $M$,
\begin{enumerate}[label=($\mathcal{A}$\arabic*)]
\setcounter{enumi}{4}
\item  one of the following holds almost surely.
\label{cond:assumption variance liminf general}
\begin{enumerate}[label=($\mathcal{A}$\arabic{enumi}-\arabic*)]
   \item $0< \liminf_{n\rightarrow \infty} \lambda_{\min} \left( \Sigma_{K,n}\right)/ |\D_n|$. \label{assumption variance liminf1}
   \item $0< \liminf_{K\rightarrow \infty} \lambda_{\min} \left(
       \Sigma_{K,n}\right)/K $.		 \label{assumption variance
       liminf2}
	\item $0< \liminf_{K, n\rightarrow \infty} \lambda_{\min} \left(
        \Sigma_{K,n}\right)/(K |\D_n|)$.\label{assumption variance liminf3}
\end{enumerate}
\end{enumerate}

The conditional centering condition \ref{cond:condcentered} is
  trivially satisfied when the theorem is applied in the context of
  conditionally unbiased estimating functions for space-time
  processes. The localization condition \ref{cond:localizationgeneral}
  in general holds trivially for the variables $\Ekn{k,n}(\l)$ while
  the measurability condition for $\EE_{\H_{k-1}}[
  Y_{n}^{(\mathcal{T}, \mathcal{I})} ]$ will typically follow from a space
  or time Markov property of the stochastic process considered. The
possibly most controversial condition is conditional
    $\alpha$-mixing \ref{cond:conditionalmdependent}. This may be quite easy to
  verify when considering stochastic processes specified in terms of
  conditional distributions, see for example
  \rwrev{Sections~\ref{sec:var} and \ref{sec:iteratedpointprocess}}. We elaborate further in
  Section~\ref{sec:alphamixing} below. The conditions of bounded
  moments \ref{cond:boundedmoments} and non-vanishing variance
  \ref{cond:assumption variance liminf general} are typical for
  central limit theorems. The more difficult condition among these
  is \ref{cond:assumption variance liminf general} which is often simply
  left as an assumption. However, we are able to give a practically
  verifiable sufficient condition for \ref{cond:assumption variance
    liminf general} in Section~\ref{sec:verifyspacetimecar}.

\subsubsection{Criteria for $\alpha$-mixing}\label{sec:alphamixing}

The condition of $\alpha$-mixing is translated into bounds for
covariances using inequalities originally due to \cite{davydov:68} and later on
\cite{rio:93} and \cite{doukhan:94} where we use the conditional
version provided in \cite{yuan2013some}. \rwrev{It is possible to
  alternatively use $\beta$-mixing to bound covariances
\citep{heinrich2013absolute} or to impose decay
conditions directly on covariances as in \cite{doukhan:lang:16}. We nevertheless} find that $\alpha$-mixing is a well-established, concise and intuitively
appealing measure of spatial dependence. Specific examples of
$\alpha$-mixing processes are $m$-dependent processes, Gaussian
processes with sufficiently fast decay of the correlation function
\citep{doukhan:94}, Neyman-Scott processes with sufficiently fast
decay of the cluster density \citep{waagepetersen:guan:09} and Cox
point processes with $\alpha$-mixing random intensity
function. Recently \cite{poinas:delyon:lavancier:19} further established
$\alpha$-mixing for associated point processes which include the
important class of determinantal point processes. In our setting, for example, for the
iterated point process in Section~\ref{sec:iteratedpointprocess}, we
may relax the assumption that the $B_k$ are Poisson conditional on the
past by the assumption that the $B_k$ are conditionally
$\alpha$-mixing Neyman-Scott, Cox, or determinantal point
processes.

\subsection{Statement of theorem}\label{s:cltstate}

Our central limit theorem is
\begin{theorem}\label{thm:clt}
Let $(T_{K,n})_{K,n \in \N}$ be a sequence of $q$-dimensional
statistics of the form \eqref{eq:TWsum} and assume that condition~\ref{cond:condcentered} holds. Assume
further that either
\begin{enumerate}[label=(\roman*)]
\item \label{casespace} conditions~\ref{cond:spacelocalized}, \ref{cond:spatial conditional independence},
  \ref{cond:boundedmoment1}, and \ref{assumption variance liminf1}
  hold, $K$ is fixed, $n \to \infty$ and $\lim_{n \to \infty}|\D_n|=\infty$,
\item \label{casetime} conditions~\ref{cond:timelocalized}, \ref{cond:temporal conditional independence},
  \ref{cond:boundedmoment2}, and \ref{assumption variance liminf2}
  hold, $n$ is fixed, and $K\to \infty$, or
\item \label{caseboth}
  conditions~\ref{cond:spacetimelocalized}, \ref{cond:space-time
    conditional independence}, \ref{cond:boundedmoment3}, and
  \ref{assumption variance liminf3}  hold, $n,K \to \infty$,
  $\lim_{n \to \infty}|\D_n|= \infty$, and $K=O(|\D_n|^C)$ where $0<C< \min\left(\tau, \frac{1}{d} -
  2\tau\right)$ for a $0< \tau <\min\{1/(4d),-1/(2 \eta)\}$.
\end{enumerate}
Then for any $t \in \R^q$,
\begin{align*}
\EE_{\H_{0}} \left [\exp( i t^\T \Sigma_{K,n}^{-\frac{1}{2}} T_{K,n} ) \right ] \rightarrow \exp(\|t\|^2/2)
\end{align*}
almost surely.
\end{theorem}
Examples of applications of~\autoref{thm:clt} are provided in
  Section~\ref{sec:applications}  while a proof of the theorem is
  given in Section~\ref{sec:proofofclt}. We conclude this section with a few remarks.

\begin{remark}\label{rem:spacetimetempering}
When $K$ increases in case~\ref{caseboth}, dependence  may
  propagate over space which may dilute the effect of increasing spatial domain. This
  is particularly reflected by the localization condition
  \ref{cond:spacetimelocalized}. Therefore,
  to exploit spatial mixing in case \ref{caseboth} we need to temper  increase of time $K$
  relative to increase of space $|\D_n|$.
\end{remark}
\begin{remark}
The conclusion of~\autoref{thm:clt} is the same for all cases
  \ref{casespace}-\ref{caseboth} and in practice we may not need the
  more complex setting of \ref{caseboth} with $n$ and $K$ jointly tending to
  infinity. The case of $K$ fixed, $n \rightarrow \infty$, is for example relevant in rain forest ecology where spatially large data sets are collected at a modest number of census times.
\end{remark}

\begin{remark}\label{rem:conditioning}
There exists some measurable function $X_0$ so that $\H_0= \sigma(X_0)$.
Write $\phi_{K,n}(x_0)$ for $\EE [ \exp(i t^\T  \Sigma_{K,n}^{-\frac{1}{2}} T_{K,n}) |X_0=x_0]$ where expectation is with respect to a regular conditional
distribution of the $\Ekn{k,n}$, ${1 \le} k \le K$, $n \in \N$, given $X_0=x_0$. Then $\EE_{\H_{0}}[ \exp(i t^\T  \Sigma_{K,n}^{-\frac{1}{2}} T_{K,n}) ]= \phi_{K,n}(X_0)$ almost
surely. In our proof we thus show that $\phi_{K,n}(x_0)$ converges to
$\exp(\|t\|^2/2)$ for {almost} all values $x_0$ of $X_0$.
\end{remark}

\begin{remark}
The result in case (ii) ($K \rightarrow \infty$, $n$ fixed) could
alternatively under the same conditions be established using existing martingale central limit
theorems \cite[e.g.][]{brown1971martingale,hall:heyde:14}. To apply
these theorems it is required to control a sum of conditional
variances. This, however, involves all the results from
Sections~\ref{sec:boundsspatialcovariances}-\ref{sec:bounds},
Section~3 
 in the supplementary material, and derivations similar to those used
to handle the term $A_1$ in our proof. We therefore find it more
natural to use our unified framework for all cases \ref{casespace}-\ref{caseboth}.
\end{remark}

\begin{remark}
In the proof of \autoref{thm:clt}, we need to bound a large
number, of the order $K^4|\D_n|^4$, of certain covariances, see
Section~\ref{sec:proof} and the supplementary material. The
conditional centering property \ref{cond:condcentered} is crucial in
this regard since it immediately reduces the number of the non-zero
covariances to be of the order $K^3|\D_n|^4$. However, it does not
affect the order $|\D_n|^4$. This is why spatial mixing assumptions
are required in case \ref{casespace} and \ref{caseboth}.
 \end{remark}

 \section{Verification of conditions for application examples}\label{sec:verification}

In this section we discuss verification of central limit theorem conditions for the application examples of Section~\ref{sec:applications}.

\subsection{\rwrev{A vector autoregression of order one}}\label{sec:toyverification}
\rwrev{
For the vector autoregression we already discussed the conditional centering and
localization conditions
in Section~\ref{sec:simpleintrocond}. We now further discuss the
conditions \ref{cond:spatial conditional independence},
  \ref{cond:boundedmoment1}, and \ref{assumption variance liminf1} for
  the fixed time increasing space case \ref{casespace}. The condition \ref{cond:spatial conditional independence} regarding decay of
$\alpha$-mixing coefficients may look a bit
complicated because it involves conditioning on the initial
configuration $\H_0 =\sigma(X_0)$. However, for the vector autoregression, the spatial $\alpha$-mixing coefficients for the $\nu_k$ are the same for any
time step $k$ and do not depend on the initial condition. Hence we
just need polynomial decay of the common $\alpha$-mixing
coefficients. Each $E_{k,n}(\l)$ is a
finite sum involving components of the $\nu_k$,
$k=0,\ldots,K$. Therefore for fixed time horizon $K$, \ref{cond:boundedmoment1} holds provided
the $\nu_l(\l)$ have moments of
sufficiently high order (e.g.\ if they are Gaussian).

Regarding the variance condition \ref{assumption variance liminf1},
the conditional variance of $T_{K,n}$ given $\H_0$ is bounded below by \[ \Var_{\H_0}[\sum_{\l=1}^n \bar X_{0,n}(\l) \nu_0(\l)] = \sum_{\l=1}^n\sum_{\j=1}^n c(\l-\j) \bar X_{0,n}(\l) \bar X_{0,n}(\j) = \sum_{h=0}^{n-1} c(h) \sum_{\l,\j:|\l-\j|=h} \bar X_{0,n}(\l) \bar X_{0,n}(\j),\]
where, assuming the $\nu_l(\l)$ have finite second moments,
$c(\l-\j)=\Cov[\nu_1(\l),\nu_1(\j)]$.
Without further assumptions on
the initial condition it not possible to conclude whether
\ref{assumption variance liminf1} is true or not. If we assume $X_0$
is distributed as the $\nu_k$ and the covariance is positive,
$c(\cdot)>0$, then 
then the variance can further be bounded below by $n c(0) \Var \bar X_{0,n}(1)$
 so that  $\lim_{n \rightarrow \infty}
n^{-1}\Var[\sum_{\l=1}^n \bar X_{0,n}(\l) \nu_0(\l)]$ is positive.

For the vector autoregression example it seems easier to reason about the
unconditional variance of $T_{K,n}$ than the conditional variance
given $\H_0$. On the other hand, strong assumptions on $X_0$ are
needed, and from a practical point of view it may be difficult to
estimate the unconditional variance when the distribution of $X_0$ is
not known. One may object that using our central limit theorem,
statistical inferences are conditional on $X_0$. But the same can be
said about inferences in regression models that are conditional on the observed covariates.

In case of increasing time fixed space (fixed $n$) we may consider the
vector form $X_{k,n}=aB_n X_{k-1,n}+\nu_{k,n}$ of the autoregression where $B_n$ is an
$n\times n$ matrix, $X_{k,n}=(X_{k,n}(\l))_{\l=1}^n$ and similar for $X_{k-1,n}$ and $\nu_{k,n}$.
By the Perron-Frobenius theorem, $aB_n$ has maximal absolute eigen value
$|a|$. Assuming $|a|<1$ allows to control covariances and
expectations over time and to use derivations similar to those
detailed in Section~\ref{sec:verifyspacetimecar}. We hence omit
further considerations for the increasing time fixed space case here.}

 \subsection{Discrete time spatial birth-death point process}\label{sec:verifyiteratedpointprocess}

 Continuing from Section~\ref{sec:applications}, for any $A \subseteq \R^d$ and $l \le k \in \Z_{+}$, we define
 $\H_{l,k,A}=  \sigma( (X_{k'} \cap A),l\leq k' \leq k)$. We consider first case \ref{casespace} of \autoref{thm:clt} and assume that the $X_k$ are observed within an increasing sequence of observation windows $W_n$. Accordingly we introduce the subscript $n$ for the index set $\D$, the variables $\Ekn{k}(\l)$, and the score function $T_k^W$.

By the first Bartlett identity (unbiasedness of score
  functions), $E_{k,n}(\l)$ is conditionally centered. Since $B_k\cup I_k=X_{k} \setminus X_{k-1}$ and $S_k=X_{k} \cap
X_{k-1}$, $E_{k,n}(\l)$ is
$\H_{0,k,c(\l)}$ measurable. Defining $Y_{n}^{(\mathcal{T},
  \mathcal{I})}$ as in \ref{cond:localizationgeneral},  $\EE_{\H_{k-1}}[ Y_{n}^{(\mathcal{T}, \mathcal{I})}]$
is a function only of $X_{k-1}$  restricted to $(\cup_{\i \in
  \mathcal{I}}c(\i))\oplus (K-k)\omega$ and is hence $\H_{0,k-1,(\cup_{\i \in
  \mathcal{I}}c(\i))\oplus (K-k)\omega}$ measurable, which means that
\ref{cond:spacelocalized} holds (here $A \oplus s=\{x \in \R^d|
\dd(\{x\},A) \le s\}$ for $A \subset \R^d$ and $s>0$).
Condition~\ref{cond:spatial conditional independence} holds due to conditional independence of the Bernouilli variables $R_u$ and since by independence properties of Poisson processes, $(B_k \cap I_k \cap
A)$ and $(B_k \cap I_k \cap B)$ are independent
given $\H_{k-1}$, $k=1,\ldots,K$, whenever $\dd(A,B) > 0$. It is easy to check that \ref{cond:boundedmoment1} holds and assuming further that \ref{assumption variance liminf1} holds, case \ref{casespace} of our \autoref{thm:clt} is applicable to
\[
T_{K,n}= \sum_{k=1}^K T_{k}^{W_n} = \sum_{k=1}^K \sum_{\l \in \D_n} \Ekn{k,n}(\l),
\]
for any fixed $K$. The validity of \ref{assumption variance liminf1}
depends on $X_0$ that e.g.\ must not be confined to a bounded region.

Imposing further conditions it is also possible to use the joint
space-time asymptotic case \ref{caseboth} of
\autoref{thm:clt}. Suppose  that the immigrant intensity $\rho$ is
strictly positive, ensuring that the `forest' never becomes
extinct. Also for $u$ in $X_0$ or $u \in \cup_{k \ge 1} I_k$ let $C_u$ denote the cluster consisting of $u$ and all its descendants, c.f.\ the birth-death dynamics for $X_k$, $k\ge 0$. If the probability of survival is small enough compared to the birth intensities, $C_u$ will be finite almost surely and hence stochastically bounded in space and time. Thereby techniques similar to those in \cite{cheysson2022spectral} can be used to establish the required $\alpha$-mixing for case \ref{caseboth} of \autoref{thm:clt}.

\subsection{\rwrev{Graphical} autoregressive model}\label{sec:verifyspacetimecar}

For this example we consider a fixed spatial domain $\D$ and hence for
notational simplicity skip the index $n$ which is only pertinent for
increasing spatial domain asymptotics.

We define \rwrev{$\H_{k, l, A}=\sigma(\{ X_{k'}(\l): \l \in A \cap \D, k\leq k' \leq l \})$
for $k \le l\in\Z_{+}$ and  $A \cap \D \neq \emptyset$. If $A \cap
\D=\emptyset$ we define $\H_{k, l, A}$ to be the trivial $\sigma$-algebra.} The conditioning
in \eqref{eq:condspec} is then equivalent to conditioning on the
history $\H_{k-1}=\H_{0,k-1,\D}$ and the neighborhood information
$\H_{k,k,\D\backslash\{\l\}}$. Given $\H_{k-1}$, by Brook's factorization Lemma \citep{besag1974spatial}, $X_{k}$ is a Gaussian Markov random field with mean vector $\mu_{k-1} = \beta \xi_{k-1}^{\mathrm{temp}} $ and precision matrix
$
  Q = A \left(I_{|\D|} - \gamma B_{0} \right),
$
where $I_{|\D|}$ is the $|\D|\times |\D|$  identity matrix. We assume $Q$ is positive definite and symmetric.

Given $\H_{k-1}$,
\begin{equation}\label{eq:xkgmrf}
  X_{k} \sim  \beta \xi_{k-1}^{\mathrm{temp}} + Q^{-1/2} Z_{k}, \quad k\in\N,
\end{equation}
where the $Z_{k}$, $k \ge 1$, are independent standard  normal vectors.

Define for $m \in \N$ and $k \le l$ the maximal conditional
correlation coefficient
\[
  \rho_{\H_{k-1}}(\H_{0,k,\D}, \H_{(l-r)^{+},\rwrev{l-1},\D}) = \sup \big| \Corr_{\H_{{k-1}}}[Z, Y]  \big|,
\]
where the supremum is over all random variables $Z$  and $Y$ that are
respectively $\H_{0,k,\D}$ and $\H_{(l-r)^{+},\rwrev{l-1},\D}$ measurable, and
where $\EE_{\H_{k-1}}[Z^2]<\infty$ and
$\EE_{\H_{k-1}}[Y^2]<\infty$. It is shown in
  Section~5.1 
  in the supplementary material  that for a $\psi<-1$ and any $\l,\j\in\D$, $\Cov_{\H_{k-1}} \left[ X_{k}(\l), X_{\rwrev{k+m}}(\j) \right] =
O(m^{\psi})$
as $m\to \infty$ so that
\[
  \sup_{\substack{l,k\in\Z_{+} \\ l-k\geq m}}  \rho_{\H_{k-1}}
    (\H_{0,k, \D},\H_{(l-r)^{+},\rwrev{l-1}, \D}) \leq \text{ constant } m^{\psi}.
\]
Since \cite[see inequality~(1.12) in][]{bradley2005basic}
\[
  \alpha_{\H_{k-1}}(\H_{0,k,\D}, \H_{(l-r)^{+}, \rwrev{l-1},\D}) \leq \frac{1}{4} \rho_{\H_{k-1}}(\H_{0,k,\D}, \H_{(l-r)^{+},\rwrev{l-1},\D}),
\]
we have
\begin{align*}
  \rwrev{\alpha_r}(m) \leq \frac{1}{4} \sup_{\substack{l,k\in\Z_{+} \\ l-k\geq m}}   \rho_{\H_{k-1}} (\H_{0,k, \D}, \H_{(l-r)^{+},\rwrev{l-1}, \D}) \leq \text{ constant } m^{\psi},
\end{align*}
which means that $\EE_{\H_{0}}\left[ \rwrev{\alpha_r}(m) \right] ^{\frac{p-2}{p}}$
is $O(m^{\psi})$ and hence condition~\ref{cond:temporal conditional independence} holds.

By definition, $\xi_{k-1,n}^{\mathrm{temp}}(\l)$ and $ \xi_{k,n}^{\mathrm{spat}}(\l)$ are  measurable with respect to $\H_{(k-r)^{+},k-1,\D}$ and $\H_{k,k,\D_{n}\backslash\{\l\}}$.
The conditional centering~\ref{cond:condcentered} holds because
\[
  \EE_{\H_{k-1}}[ E_{k,n}(\l) ] = \EE_{\H_{k-1}}\left[ \EE_{\sigma(\H_{k-1}\cup \H_{k,k,\D_{n}\backslash\{\l\}})}\left[ \epsilon_{k,n}(\l) \right] \left(\xi_{k-1,n}^{\mathrm{temp}}(\l), \xi_{k,n}^{\mathrm{spat}}(\l)\right)^{\T} \right]=0.
\]
Hence, $T_{K,n}$ provides an unbiased estimating function.
In addition, $E_{k,n}(\l)$ is $\H_{(k-r)^{+},k,\D}$ measurable and for
any $\mathcal{T}\subseteq \{k,k=1,\ldots, K\}$, $\mathcal{I}\subset \D$, and measurable function $f:\R^{2|\mathcal{T}||\mathcal{I}|} \to \R$, $\EE_{\H_{k-1}}[f((E_{l,n}(\l))_{l\in\mathcal{T},\l\in\mathcal{I}})]$ is $\H_{(k-r)^{+}, k-1, \D}$ measurable. This means that the localization condition~\ref{cond:timelocalized} holds.

It is shown in Section~5.3 
 of the supplementary material that
\[
  \Sigma_{K,n} = \Var_{\H_{0}}[T_{K,n}] = \mathrm{diag}(\lambda^{(1)}_{K,n}, \lambda^{(2)}_{K,n}),
\]
where $\lambda_{K,n}^{(1)}\geq K \mathrm{trace} \left( B_{1,n}^{2} \right)$ and $\lambda_{K,n}^{(2)} = 2 K  \mathrm{trace}( B_{0,n}^{2})$.
We assume $\mathrm{trace}( B_{0,n}^{2} )>0$ and
  $\mathrm{trace}( B_{1,n}^{2} )>0$.
Then, $\lim\inf_{K\to\infty}  \lambda^{(1)}_{K,n} / K > 0$ and  $\lim\inf_{K\to\infty}  \lambda^{(2)}_{K,n} / K > 0$ and hence condition~\ref{assumption variance liminf2} holds.

Following Section~5.4 
 in the supplementary
material, condition~\ref{cond:boundedmoment2} holds for $\nu=6+\epsilon$ and $\epsilon\geq 3(p-2)$.
Therefore, the case~\ref{casetime} of \autoref{thm:clt} is applicable for the statistic $T_{K,n}$.

\section{Proof of the central limit theorem}
\label{sec:proofofclt}

The following
  Sections~\ref{sec:boundsspatialcovariances}-\ref{sec:bounds} contain
  prerequisites needed for the proof in Section~\ref{sec:proof}.

\subsection{Mixing covariance inequalities}\label{sec:boundsspatialcovariances}

Let $Y$ and $Z$ be random variables with
$\EE_{\H_{k}}[|Y|^{p}]<\infty$ a.s.\ and $\EE_{\H_{k}}[|Z|^{p}]<\infty$
a.s.\ for some $p>2$ and $k \in \Z_+$. Then \cite[][Theorem~3.4]{yuan2013some},
\begin{equation}\label{eq:covalphabound}
  \big| \Cov_{\H_{k}} \left[ Y, Z  \right] \big| \leq 8 \left( \EE_{\H_{k}}[|Y|^{p} ] \EE_{\H_{k}}[|Z|^{p}] \right)^{\frac{1}{p}}
  \alpha_{\H_{k}}\left(\sigma(Y), \sigma(Z) \right)^{\frac{p-2}{p}},
\end{equation}
where $\sigma(Y),\sigma(Z) \subseteq\mathcal{F}$ denote the $\sigma$-algebras
generated by $Y$ and $Z$.
By an application of
H\"{o}lder's inequality,
\begin{align}
  &\big| \EE_{\H_{0}}  \left[ \Cov_{\H_{k}} \left[ Y, Z  \right] \right] \big| \leq  
    8 \EE_{\H_{0}}\Big[ \left( \EE_{\H_{k}}[|Y|^{p} ] \EE_{\H_{k}}[|Z|^{p}] \right)^{\frac{1}{p}} \alpha_{\H_{k}}\left(\sigma(Y), \sigma(Z) \right) ^{\frac{p-2}{p}}  \Big] \le  \nonumber\\
                       &8  \EE_{\H_{0}}\Big[ \left(\EE_{\H_{k}}\left[|Y|^{p}\right]^{\frac{1}{p}} \right)^{p} \Big]^{\frac{1}{p}}\EE_{\H_{0}}\Big[ \left(\EE_{\H_{k}}\left[|Z|^{p}\right]^{\frac{1}{p}} \right)^{p} \Big]^{\frac{1}{p}}
  \EE_{\H_{0}}\left[ \left( \alpha_{\H_{k}}\left(\sigma(Y), \sigma(Z) \right) ^{\frac{p-2}{p}} \right)^{\frac{p}{p-2}} \right]^{\frac{p-2}{p}} \nonumber\\
  &= 8\EE_{\H_0}\left[|Y|^{p}\right]^{\frac{1}{p}}\EE_{\H_0}\left[|Z|^{p}\right]^{\frac{1}{p}}\EE_{\H_{0}} \left[\alpha_{\H_{k}}\left(\sigma(Y), \sigma(Z) \right)  \right] ^{\frac{p-2}{p}}. \label{eq:covalphabound3}
\end{align}

\subsection{Bounds on covariances}\label{sec:bounds}

In the proof of~\autoref{thm:clt} in Section~\ref{sec:proof},
we need to bound covariances of the form
\begin{align} \label{eq:covarianceofA1terms}
\Cov_{\H_{0}} & \left[ E_{k,n}(\l) E_{l,n}(\j), E_{k',n}(\l') E_{l',n}(\j') \right]
\end{align}
for $1\leq k,k'\leq K$ and $k\leq l\leq K, k'\leq l'\leq K$ in the
one-dimensional case $q=1$.

The covariance~\eqref{eq:covarianceofA1terms}
 can be rewritten as
\begin{align*}
 &\EE_{\H_{0}} \left[ \Cov_{\H_{k^{*}-1}}\left[ E_{k,n}(\l) E_{l,n}(\j), E_{k',n}(\l') E_{l',n}\j')  \right] \right]  \\
   &+  \Cov_{\H_{0}}\left[ \EE_{\H_{k^{*}-1}} [ E_{k,n}(\l) E_{l,n}(\j) ], \EE_{\H_{k^*-1}}[ E_{k',n}(\l') E_{l',n}(\j') ]  \right],
\end{align*}
where $k^{*} = \max\{k,l,k',l'\}$.
If precisely one index coincides  with $k^*$,  then by
  conditional centering, the covariance~\eqref{eq:covarianceofA1terms} is
zero. Thus, the covariance~\eqref{eq:covarianceofA1terms} can only be
non-zero in the cases
a) $k^*=k=l>l'\geq k'$ or $k^*=k'=l'>l \geq k$,
b) $k^*=l=l'>k,k'$,
c) $k^*=k=l=l'>k'$ or $k^*=l'=k'=l>k$, and
d) $k^*=k=l=k'=l'$.
To derive bounds in the cases a)-d) we exploit recursively the
space-time structure, conditional centering, and conditional
$\alpha$-mixing. We  consider below  a) for cases~\ref{casespace}-\ref{caseboth}. The remaining simpler b)-d) are covered in Section~3 in the supplementary material.
The obtained bounds for \eqref{eq:covarianceofA1terms} are summarized
  in~\autoref{tab:covarboundssummary}.
\begin{table}
\centering
\caption{Bounds on the covariance~\eqref{eq:covarianceofA1terms} in
  different cases.}
\label{tab:covarboundssummary}
\begin{tabular}{ccc}
\hline
a) & case~\ref{casespace} & $8 M_{\epsilon}^{4} K\bar{\alpha}(\dd(\{\l, \j\}, \{\l', \j'\}) - KR/2)$\\
 & case~\ref{casetime} & $8 M_{\epsilon}^{4}\left( \tilde
                         \alpha\left( k -1 - l'
                         \right)+I[k'=l']\sum_{i=1}^{l'-1}\tilde
                         \alpha\left( k -1 - i \right)\right)$ \\
& case~\ref{caseboth} & $8 M_{\epsilon}^{4}  \bigg ( \tilde{\alpha}(k-1-l'; d(\{\l,\j\}, \{\l', \j'\})-R)$\\
& & $  +
    I[k'=l'] \sum_{i=1}^{l'-1}\tilde{\alpha}(k-1-i; d(\{\l,\j\}, \{\l',\j'\})-(l'+1-i)R/2) \bigg )$ \\
b) & case~\ref{casespace} + \ref{caseboth} & $8 M_{\epsilon}^{4} \bar{\alpha}(\dd(\{\l, \j\}, \{\l', \j'\}) - R)$ \\
 & case~\ref{casetime} & $ 2M_{\epsilon}^{4} $ \\
c) & case~\ref{casespace} + \ref{caseboth} & $8 M_{\epsilon}^{4} \bar{\alpha}(\dd(\{\l, \j\}, \{\l', \j'\}) - R)$ \\
 & case~\ref{casetime} & $ 2M_{\epsilon}^{4} $ \\
d) & case~\ref{casespace} + \ref{caseboth} & $8 M_{\epsilon}^{4} K\bar{\alpha}(\dd(\{\l, \j\}, \{\l', \j'\}) - KR/2)$ \\
 & case~\ref{casetime} & $ 2 M_{\epsilon}^{4} $ \\\hline
\end{tabular}
\end{table}

\subsubsection{Bounds in case a) \ref{casespace}-\ref{caseboth}}
Consider the case $k^*=k=l>l' \geq k'$ ($k^*=k'=l'>l \geq k$ is handled the same way). Then
$
  \Cov_{\H_{k-1}}\left[ E_{k,n}(\l) E_{k,n}(\j), E_{k',n}(\l') E_{l',n}(\j')  \right] = 0
$
and hence the covariance~\eqref{eq:covarianceofA1terms} reduces to
\begin{equation}\label{eq:reducedcov}
  \Cov_{\H_{0}} \left[ \EE_{\H_{k-1}}\left[ E_{k,n}(\l) E_{k,n}(\j) \right], E_{k',n}(\l') E_{l',n}(\j') \right].
\end{equation}
We rewrite \eqref{eq:reducedcov} as
\begin{align}
&\EE_{\H_{0}} \left[ \EE_{\H_{k-1}}\left[ E_{k,n}(\l) E_{k,n}(\j)
                \right]\left(E_{k',n}(\l')
                E_{l',n}(\j')-\EE_{\H_{0}}[E_{k',n}(\l')
                E_{l',n}(\j')]\right)\right] \nonumber \\
= &\EE_{\H_{0}} \left[ \EE_{\H_{k-1}}\left[ E_{k,n}(\l) E_{k,n}(\j)
    \right]\sum_{i=1}^{l'}\left(\EE_{\H_{i}}[E_{k',n}(\l')
    E_{l',n}(\j')]-\EE_{\H_{i-1}}[E_{k',n}(\l')
    E_{l',n}(\j')]\right)\right] \nonumber \\
= &\sum_{i=1}^{l'}\EE_{\H_{0}} \left[
    \Cov_{\H_{i-1}}\left(\EE_{\H_{k-1}}[ E_{k,n}(\l) E_{k,n}(\j)],
    \EE_{\H_{i}}[E_{k',n}(\l') E_{l',n}(\j')]\right)\right]. \label{eq:termi}
\end{align}
If $k'<l'$ then $\EE_{\H_{l'-1}}[E_{k',n}(\l') E_{l',n}(\j')]=0$ and therefore $\EE_{\H_{i}}[E_{k',n}(\l') E_{l',n}(\j')]=0$ for all $i\leq l'-1$ meaning that the only non-zero term in the previous expression is the one for $i=l'$, i.e.\
\begin{equation}\label{eq:il'} \EE_{\H_{0}} \left[
    \Cov_{\H_{l'-1}}\left(\EE_{\H_{k-1}}[ E_{k,n}(\l) E_{k,n}(\j)],
      E_{k',n}(\l') E_{l',n}(\j')\right) \right].\end{equation}

Considering \ref{caseboth}, if condition~\ref{cond:spacetimelocalized} holds, $\EE_{\H_{k-1}} [ E_{k,n}(\l) E_{k,n}(\j)]$ and $E_{k',n}(\l') E_{l',n}(\j')$ are measurable with respect to respectively \apMOD{$\H_{(k-r)^{+},k-1,C(\l)\cup C(\j)}$ and $\H_{0,l',C(\l')\cup C(\j'))}$}. Since
\begin{align*}
& \alpha_{\H_{l'-1}}  \left(\H_{0,l',C(\l')\cup C(\j')},\H_{(k-r)^{+},k-1,C(\l)\cup C(\j)}\right)\\
\leq & \alpha_{\rwrev{r}; k,2R^2,2R^2}(k-1-l'; d(C(\l)\cup C(\j), C(\l')\cup C(\j')))\\
\leq & \alpha_{\rwrev{r}; k,2R^2,2R^2}(k-1-l'; d(\{\l,\j\}, \{\l', \j'\})-R),
\end{align*}
we can use \eqref{eq:covalphabound3} to bound \eqref{eq:il'} and hence \eqref{eq:reducedcov} by
\begin{multline*}8\EE_{\H_0} \left [ \left |\EE_{\H_{k-1}}[ E_{k,n}(\l) E_{k,n}(\j)]\right|^p \right]^{\frac{1}{p}}\EE_{\H_0} \left [ |E_{k',n}(\l') E_{l',n}(\j')|^p \right]^{\frac{1}{p}}\\
\times\EE_{\H_0} \left [\alpha_{r; k,2R^2,2R^2}(k-1-l'; d(\{\l,\j\}, \{\l', \j'\})-R)\right]^{\frac{p-2}{p}}.
\end{multline*}
By Jensen's inequality, we can bound $\EE_{\H_0}\! \left [ \left
    |\EE_{\H_{k-1}}[ E_{k,n}(\l) E_{k,n}(\j)]\right|^p\right]^{\frac{1}{p}}$ by $\EE_{\H_0} \! \left[\left | E_{k,n}(\l) E_{k,n}(\j)]\right|^p \right]^{\frac{1}{p}}$ and then by $M_{\epsilon}^2$ using equation (2) 
in the supplementary material. Using the same reasoning on the second factor in the previous expression and using condition~\ref{cond:space-time conditional independence} on the third factor yields the bound
\begin{equation}\label{eq:boundil'}  8 M_{\epsilon}^4 \tilde \alpha(k-1-l',d(\{\l,\j\}, \{\l', \j'\})-R). \end{equation}
for the covariance \eqref{eq:reducedcov}.

Consider next $l'=k'$. In this case we distinguish between $i=l'$ and
$i<l'$ for the terms in \eqref{eq:termi}. For the term with $i=l'$ we
by similar arguments as above again obtain the bound
\eqref{eq:boundil'}. For $i<l'$, if \ref{cond:spacetimelocalized}
holds, $\EE_{\H_{i}}[E_{l',n}(\l')E_{l',n}(\j')]$ is measurable with respect to
\apMOD{$\H_{(i+1-r)^{+},i,((l'-i)C(\l')\cup C(\j'))}$}. Thus, since
\begin{align*}
& \alpha_{\H_{i-1}}  \left(\H_{(k-r)^{+},k-1,C(\l)\cup C(\j)}, \H_{(i+1-r)^{+},i,(l'-i)(C(\l')\cup C(\j'))}\right)\\
\leq & \apMOD{\alpha_{r; k,2R^2,2((k-i)R)^2}}(k-1-i; d(C(\l)\cup C(\j), (l'-i)(C(\l')\cup C(\j'))))\\
\leq & \apMOD{\alpha_{r; k,2R^2,2((k-i)R)^2}}(k-1-i; d(\{\l,\j\}, \{\l', \j'\})-(l'-i+1)R/2)
\end{align*}
from (2) 
 in the
supplementary material, \eqref{eq:covalphabound3}, and condition~\ref{cond:space-time conditional independence}
we obtain the bound $8 M_{\epsilon}^4 \tilde \alpha(k-1-i),d(\{\l,\j\}, \{\l', \j'\})-(l'+1-i)R/2)$.

Collecting everything, \eqref{eq:reducedcov} is bounded in case (iii) by
\begin{multline*}
8 M_{\epsilon}^{4} \Big ( \tilde \alpha(k-1-l'; d(\{\l,\j\}, \{\l',
  \j'\})-R)\\  +I[k'=l']\sum_{i=1}^{l'-1}\tilde \alpha(k-1-i; d(\{\l,\j\}, \{\l', \j'\})-(l'+1-i)R/2) \Big ).
\end{multline*}
In case (i) with $K$ fixed we can use $l'\le K$ and simplify the bound to $8 M_{\epsilon}^{4} \bar{\alpha}(d(\{\l,\j\}, \{\l', \j'\})-KR/2)$
  and in case (ii) with $n$ fixed we can simplify to $8 M_{\epsilon}^{4} \left ( \tilde \alpha(k-1-l')+I[k'=l']\sum_{i=1}^{l'-1}\tilde \alpha(k-1-i) \right )$.

\subsection{Proof of \autoref{thm:clt}}\label{sec:proof}


\begin{proof}
By an extension of the Cram{\' e}r-Wold  device
\cite[Lemma~2.1 in][]{biscio:poinas:waagepetersen:17}, it suffices
to verify Theorem~\ref{thm:clt} in the univariate case  $q=1$. Then $\Sigma_{K,n}=\sigma_{K,n}^2 = \Var_{\H_{0}}[ T_{K,n} ]$ is scalar and by conditional centering~\ref{cond:condcentered},
\begin{align*}
\sigma_{K,n}^2 =& \sum_{k,k'=1}^K \sum_{\l,\l' \in \D_n} \EE_{\H_{0}}[ E_{k,n}(\l) E_{k',n}(\l') ]
  = \sum_{k=1}^K \sum_{\l,\l' \in \D_n} \EE_{\H_{0}}[ E_{k,n}(\l) E_{k,n}(\l') ].
\end{align*}

We introduce  spatial and temporal truncation distances $\mkn$ and
$\lkn$ as follows. In case \ref{casespace}, $\lkn=\infty$ and $\mkn=|\D_n|^\tau$ with $0 < \tau < \min\{1/(4d),-1/(2 \eta)\}$ so that as $n\to \infty$,
\begin{equation}\label{eq:mkn}
  \mkn\to\infty, \quad  \frac{\mkn^{2d}}{|\D_n|^{1/2}} \to 0 \quad \text{and} \quad \mkn^{\eta}|\D_n|^{1/2} \to 0.
\end{equation}
In case \ref{casetime}, $\mkn=\infty$ and $\lkn=K^{\tilde \tau}$ with $0<\tilde \tau<-1/(2\psi)$ so that as $K \rightarrow \infty$,
\begin{equation}\label{eq:lkn} \lkn \rightarrow \infty,  \quad  \frac{\lkn}{|K|^{1/2}} \rightarrow 0 \quad
  \text{and} \quad  \lkn^{\psi}|K|^{1/2} \rightarrow 0 .\end{equation}
In case \ref{caseboth}, we choose $\mkn=|\D_n|^\tau$ and $\lkn=K^{\tilde \tau}$ with $\tau$ and $\tilde \tau$ as above, so that as $n,K \rightarrow \infty$, both \eqref{eq:mkn} and \eqref{eq:lkn} are satisfied.

 We further define $\bar{S}_{K,n} = {{T}_{K,n}}/{V_{K,n}}$ where
\[
  V_{K,n}^2 = \sum_{k=1}^K \sum_{(\l,\j) \in \mathcal{J}_{K,n}} \EE_{\H_{0}}[ E_{k,n}(\l) E_{k,n}(\j)],
\]
and $\mathcal{J}_{K,n} = \{(\l, \j): \l, \j\in\D_n,  \dd(\l,\j) \leq
\mkn\}$. Then $V_{K,n}^2$ coincides with $\sigma_{K,n}^2$ in case
\ref{casetime} and $V_{K,n}^2$ is a spatially truncated version of $\sigma_{K,n}^2$ in
cases \ref{casespace} and \ref{caseboth}. In Section~4 
of the supplementary
material it is shown that, uniformly in $K$,
\begin{equation}\label{eq:approxvar}
  \left|1 - V_{K,n}^2/\sigma_{K,n}^2 \right| \xrightarrow[n\rightarrow \infty ]{a.s.} 0.
\end{equation}

We next need to define a quantity $\bar{S}_{K,n}(k,\l)$ (see paragraph below) as a sum over variables
$\Ekn{l,n}(\j)$ so that we can
show convergence to zero of $\EE_{\H_{0}}[ A_1 ]$, $\EE_{\H_{0}}[
A_2]$, and $\EE_{\H_{0}}[A_3]$ where
\begin{align*}
  A_1 &= it e^{it \bar{S}_{K,n}} \bigg( 1- \frac{1}{V_{K,n}}  \sum_{k=1}^K \sum_{\l \in \D_n} E_{k,n}(\l) \bar S_{K,n}(k,\l) \bigg), \\
  A_2 &= \frac{e^{it \bar{S}_{K,n}}}{V_{K,n}} \sum_{k=1}^K \sum_{\l \in \D_n} E_{k,n}(\l)  \left( 1 - it \bar{S}_{K,n}(k,\l) - \exp\left\{-it \bar{S}_{K,n}(k,\l) \right\}  \right), \\
  A_3 &= \frac{1}{V_{K,n}} \sum_{k=1}^K \sum_{\l \in \D_n} E_{k,n}(\l)  \exp\left\{it \left(\bar{S}_{K,n}- \bar{S}_{K,n}(k,\l) \right) \right\},
\end{align*}
and, following \cite{bolthausen:82},
$
  \left(it-\bar{S}_{K,n}\right) e^{it \bar{S}_{K,n}} = A_1 - A_2 -A_3.
$

If the above mentioned convergences hold then  for any $t \in \R$,
$
\EE_{\H_{0}} \left[ \left(it-\bar{S}_{K,n}\right) e^{it \bar{S}_{K,n}}
  \right ] \rightarrow 0,
$
  which by Lemma~2 in \cite{bolthausen:82} \cite[see also the discussion in][]{biscio:poinas:waagepetersen:17} implies $\EE_{\H_{0}} \left[ \exp\left(i t
    \bar{S}_{K,n} \right) \right] \rightarrow \exp(t^2/2)$ which again
implies the desired result $\EE_{\H_{0}}\left[ \exp\left( it \frac{T_{K,n}}{\sigma_{K,n}}\right) \right] \rightarrow \exp(t^2/2)$.

To handle $A_1$ we want $\EE_{\H_{0}} \left[ \sum_{k=1}^K \sum_{ \l
    \in\D_n} E_{k,n}(\l)  \bar S_{K,n}(k,\l) \right]/
V_{K,n}^2=1$ and we want a small number of terms in $\bar
S_{K,n}(k,\l)$. However, we also prefer that $S_{K,n}(k,\l)$
is close to $\bar S_{K,n}$ in order to deal with $A_3$. A suitable
compromise is given by
\begin{align*}
 \bar{S}_{K,n}(k,\l) &= \frac{1}{V_{K,n}} \sum_{(l,\j) \in \mathcal{I}_{K,n}(k,\l)} E_{l,n}(\j),
\end{align*}
where
\[
 \mathcal{I}_{K,n}(k,\l)  = \{(l, \j): k \le l \le K,  l-k \le
 \lkn, (\l,\j) \in \mathcal{J}_{K,n} \}.
\]
Compared to $\bar S_{K,n}$ we avoid all $\Ekn{l,n}(\j)$ for $l$ prior to
  $k$. Combined with the truncation by $\mkn$ or $\lkn$ in space or
  ahead in time (or both) this makes the number of terms in $\bar S_{K,n}(k,\l)$
  small. Also, omitting $\Ekn{l,n}(\j)$ terms for $l<k$ turns
    out not to be a problem for handling $A_3$.

In the remainder we consider the convergences of $\EE_{\H_{0}}[ A_1 ]$, $\EE_{\H_{0}}[
A_2]$, and $\EE_{\H_{0}}[A_3]$.

\subsubsection{Convergence of  $\EE_{\H_0}[A_1]$}

First note that
\begin{align*}
  &\EE_{\H_{0}} \left[ \frac{1}{V_{K,n}} \sum_{k=1}^K \sum_{ \l \in\D_n} E_{k,n}(\l)  \bar S_{K,n}(k,\l) \right]
  = \frac{1}{V_{K,n}^{2}} \sum_{k=1}^K \sum_{\l\in\D_n} \sum_{(l,\j) \in \mathcal{I}_{K,n}(k,\l)} \EE_{\H_{0}}[ E_{k,n}(\l) E_{l,n}(\j)]\\
= &\frac{1}{V_{K,n}^{2}}\sum_{k=1}^K \sum_{(\l,\j) \in \mathcal{J}_{K,n}} \EE_{\H_{0}}[ E_{k,n}(\l) E_{k,n}(\j)]=1,
\end{align*}
where the second  equality is due to conditional centering.
Thus,
\begin{align*}
   &\EE_{\H_{0}} \left[ \left|A_1 \right| ^2 \right]  \leq \frac{t^2}{V_{K,n}^4} \Var_{\H_{0}}\bigg[ \sum_{k=1}^K \sum_{\l \in \D_n} \sum_{(l,\j) \in \mathcal{I}_{K,n}(k,\l)} E_{k,n}(\l) E_{l,n}(\j)  \bigg] \notag\\
    =&\frac{t^2}{V_{K,n}^4} \sum_{k,k'=1}^K \sum_{\l,\l' \in \D_n} \!\!\!
       \sum_{\substack{(l,\j) \in \mathcal{I}_{K,n}(k,\l) \\ (l',\j') \in \mathcal{I}_{K,n}(k',\l')}}  \!\!\! \!\!\! \!\Cov_{\H_{0}}\left[ E_{k,n}(\l) E_{l,n}(\j), E_{k',n}(\l') E_{l',n}(\j') \right]
   = \frac{t^2}{V_{K,n}^4} \left( 2 a + b + 2 c + d \right), 
\end{align*}
where
\begin{align*}
  a &= \sum_{k=2}^{K} \sum_{l'=1}^{k-1}  \sum_{\substack{k'=1\\ l'-k'\le \lkn}}^{l'}\sum_{(\l,\j),(\l',\j')\in\mathcal{J}_{K,n}}  \Cov_{\H_{0}} \left[ E_{k,n}(\l) E_{k,n}(\j), E_{k',n}(\l') E_{l',n}(\j') \right] \\
  b &= \sum_{l=2}^{K} \sum_{\substack{k, k'= 1\\ l-k,l-k' \le \lkn}}^{l - 1}  \sum_{(\l,\j),(\l',\j')\in\mathcal{J}_{K,n}}  \Cov_{\H_{0}} \left[ E_{k,n}(\l) E_{l,n}(\j), E_{k',n}(\l') E_{l,n}(\j')  \right]\\
  c &= \sum_{k=2}^{K} \sum_{\substack{k'=1\\k-k'\le \lkn}}^{k-1}   \sum_{(\l,\j),(\l',\j')\in\mathcal{J}_{K,n}}  \Cov_{\H_{0}} \left[ E_{k,n}(\l)  E_{k,n}(\j), E_{k',n}(\l') E_{k,n}(\j')   \right]  \\
  d &= \sum_{k=1}^{K}  \sum_{(\l,\j),(\l',\j')\in\mathcal{J}_{K,n}}  \Cov_{\H_{0}} \left[ E_{k,n}(\l) E_{k,n}(\j), E_{k,n}(\l') E_{k,n}(\j') \right]
\end{align*}
are corresponding to cases a)-d) in Section~\ref{sec:bounds}.
 We control a-d in each of the cases
\ref{casespace}-\ref{caseboth} using the bounds in Table~\ref{tab:covarboundssummary}.\\
\underline{Case~\ref{casespace}}: In case \ref{casespace},
\begin{align*}
  |a| &\leq 8K^{4} M_{\epsilon}^{4}
        \sum_{(\l,\j),(\l',\j')\in\mathcal{J}_{K,n}}  \bar{\alpha}\left(\dd(\{\l,\j\},\{\l',\j'\})-KR/2\right) \\
  |b| &\leq 8K^{3} M_{\epsilon}^{4}  \sum_{(\l,\j),(\l',\j')\in\mathcal{J}_{K,n}}  \bar{\alpha}\left(\dd(\{\l,\j\},\{\l',\j'\})-R\right) \\
  |c| &\leq 8K^{2} M_{\epsilon}^{4}  \sum_{(\l,\j),(\l',\j')\in\mathcal{J}_{K,n}}  \bar{\alpha}\left(\dd(\{\l,\j\},\{\l',\j'\})-R\right) \\
  |d| &\leq 8K^{2} M_{\epsilon}^{4}  \sum_{(\l,\j),(\l',\j')\in\mathcal{J}_{K,n}}  \bar{\alpha}\left(\dd(\{\l,\j\},\{\l',\j'\})-KR/2\right).
\end{align*}
The proof now proceeds by splitting the above sums according to
whether $\dd(\l,\l') < 3m_n$ or not. We omit the details which follow
quite closely the proof of the spatial case \cite[see e.g.\
][]{biscio:waagepetersen:18} and obtain that
$a=O(K^{4} |\D_n|\mkn^{2d})$, $b=O(K^3 |\D_n|\mkn^{2d})$,
$c=O(K^2|\D_n|\mkn^{2d})$ and $d=O(K^{2}|\D_n|\mkn^{2d})$. Hence, $a,b,c,d$ are all $O(|\D_n|\mkn^{2d})$ for $K$ fixed.\\[\bsl]
%
\underline{Case~\ref{casetime}}: In case \ref{casetime},
\begin{align*}
|a|\leq & 8 M_{\epsilon}^{4}|\mathcal{D}_n|^4\sum_{k=2}^{K} \sum_{l'=1}^{k-1}  \sum_{\substack{k'=1\\ l'-k'\le \lkn}}^{l'}\left(\tilde \alpha( k -1 - l')+I[k'=l']\sum_{i=1}^{l'-1}\tilde \alpha\left( k -1 - i \right)\right) \\
\leq & 8 M_{\epsilon}^{4}|\mathcal{D}_n|^4\left( \lkn \sum_{k=2}^{K} \sum_{l'=1}^{k-1} \tilde \alpha( k -1 - l') + \sum_{k=2}^{K} \sum_{l'=1}^{k-1}\sum_{z=k-1-l'}^{k-2}\tilde \alpha(z)\right).
\end{align*}
The first term can be directly bounded by $8 M_{\epsilon}^{4}|\mathcal{D}_n|^4 \lkn K \sum_{z=0}^{\infty}\tilde \alpha(z)$. By condition~\ref{cond:temporal conditional independence}, $\tilde{\alpha}_{n}(z)=O(z^{\psi})$, $\psi< -1$, and $\sum_{z=0}^{\infty}   \tilde{\alpha}_{n}(z) \leq \text{constant} \sum_{z=0}^{\infty} z^{\psi} < \infty$.
Hence this expression is $O(K \lkn|\D_n|^{4})$.
The second term can be rewritten as
\begin{align}
& 8 M_{\epsilon}^{4}|\mathcal{D}_n|^4\sum_{z=0}^{K-2}\tilde \alpha(z)\#\{(k, l'): 1\leq l'<k\leq K, k-1-l'\leq z\leq k-2\}\nonumber\\
\leq & 8 M_{\epsilon}^{4}|\mathcal{D}_n|^4\sum_{z=0}^{K-2}\tilde
       \alpha(z)\#\{(k, t): 2\leq k\leq K, 1\leq t\leq K, t-1\leq
       z\leq k-2\} \nonumber \\
\leq & 8 M_{\epsilon}^{4}|\mathcal{D}_n|^4\sum_{z=0}^{K-2}\tilde \alpha(z)\#\{(k, t): 2+z\leq k\leq K, 1\leq t\leq z+1\}\nonumber\\
\leq & 8 M_{\epsilon}^{4}|\mathcal{D}_n|^4 K^2 \left( \frac{1}{K} \sum_{z=0}^{K-2}(z+1)\tilde \alpha(z)\right).\label{eq:a}
\end{align}
By condition~\ref{cond:temporal conditional independence},
$(z+1)\tilde
\alpha_{n}(z)\underset{z\rightarrow\infty}{\longrightarrow} 0$. Hence,
by Ces\`aro summation, $\frac{1}{K} \sum_{z=0}^{K-2}(z+1)\tilde
\alpha_{n}(z)\underset{K\rightarrow\infty}{\longrightarrow} 0$ and
\eqref{eq:a} is thus $o(|\mathcal{D}_n|^4 K^2)$. We further have $|b| \leq 2M_{\epsilon}^{4}|\D_n|^{4}l_K^2$, $|c| \leq 2M_{\epsilon}^{4}|\D_n|^{4}l_K$ and $|d| \leq 2 M_{\epsilon}^{4}|\D_n|^{4}  K$. Then $a, b, c,
d$ are all $O(K \lkn^2|\D_n|^{4})+o(K^2|\D_n|^{4})$ which is $o(K^2)$ for $n$ fixed as a consequence of \eqref{eq:lkn}.\\[\bsl]

\noindent \underline{Case~\ref{caseboth}}: In case \ref{caseboth},
\begin{align*}
|b| &\leq 8 M_\epsilon^4 \sum_{l=2}^K\sum_{\substack{k, k'= 1\\ l-k,l-k' \le \lkn}}^{l - 1}  \sum_{(\l,\j),(\l',\j')\in\mathcal{J}_{K,n}} \bar{\alpha}(d(\{\j\}, \{\j'\})-R)\\
&\leq 8 M_\epsilon^4 Kl_K^2 \sum_{(\l,\j),(\l',\j')\in\mathcal{J}_{K,n}}\bar{\alpha}(d(\{\j\},\{\j'\})-R)
\end{align*}

and

\begin{align*}
|c| &\leq 8 M_{\epsilon}^{4}\sum_{k=2}^K  \sum_{\substack{k'=1\\k-k'\le \lkn}}^{k-1}   \sum_{(\l,\j),(\l',\j')\in\mathcal{J}_{K,n}} \bar{\alpha}(d(\{\l, \j\}, \{\j'\})-R)\\
&\leq 8 M_{\epsilon}^{4}Kl_K \sum_{(\l,\j),(\l',\j')\in\mathcal{J}_{K,n}} \bar{\alpha}(d(\{\l, \j\}, \{\j'\})-R).
\end{align*}
Similarly to case \ref{casespace}, by splitting the sum according to
whether $\dd(\l,\l') < 3m_n$ or not \cite[see][]{biscio:waagepetersen:18} we get that $|b|=O(Kl_K^2|\D_n|m_n^{2d})$ and $|c|=O(Kl_K|\D_n|m_n^{2d})$.

For $|d|$ we get
$$|d| = 8 M_{\epsilon}^{4}
K^2\sum_{(\l,\j),(\l',\j')\in\mathcal{J}_{K,n}}\bar{\alpha}(\dd(\{\l,
\j\}, \{\l', \j'\}) - KR/2).$$
Bounding $|d|$ is a bit more tricky than bounding $|b|$ and $|c|$ due
to the dependence on $K$ inside the mixing coefficient. We thus give
details on how this term is handled. Following
\cite{biscio:waagepetersen:18} we split the above sum according to
whether $\dd(\l,\l') < 3m_n$ or not and conclude that there exists a constant
$\tilde C>0$ such that
\begin{align*}
&\sum_{(\l,\j),(\l',\j')\in\mathcal{J}_{K,n}}\bar{\alpha}(\dd(\{\l, \j\}, \{\l', \j'\}) - KR/2)\\
&~~~~~~~~\leq \tilde C|\D_n|m_n^{2d}\left(\sum_{z\geq 3m_n}z^{d-1}\bar{\alpha}(z - 2m_n - KR/2) + \sum_{z\leq 4m_n}z^{d-1}\bar{\alpha}(z - KR/2) \right).
\end{align*}
For the first sum, since $C<\tau$ then $K=o(|D_n|^\tau)=o(m_n)$ and we
get for sufficiently large $n$,
$$\sum_{z\geq 3m_n}z^{d-1}\bar{\alpha}(z - 2m_n - KR/2)=O\left(
\sum_{z \geq 3m_n}^{\infty}z^{\eta+d-1}\right) =o(1).$$
By splitting the second sum according to whether $z - KR/2<0$ or not we get
\begin{align*}
\sum_{z\leq 4m_n}z^{d-1}\bar{\alpha}(z - KR/2) &\leq   \!\!\!\!\!   \sum_{z\leq \lfloor  KR/2\rfloor} \!\!\!\!\! z^{d-1} + \!\!\! \!\! \!\!\!\sum_{z \leq \lfloor 4m_n-KR/2\rfloor}  \!\!\! \!\!\! \!\! (z+KR)^{d-1}\bar{\alpha}(z)\\
& \leq  \left(\frac{KR}{2}\right)^2 + \sum_{z=0}^\infty (z+KR/2)^{d-1}\bar{\alpha}(z).
\end{align*}
The term $(z+KR/2)^{d-1}$ consists of monomials of the form $z^iK^{d-1-i}$ for all $i\leq d-1$ and $\sum_{z=0}^\infty z^i\bar{\alpha}(z)$ is finite by assumption \ref{cond:space-time conditional independence}. Thus
$$\left(\frac{KR}{2}\right)^2 + \sum_{z=0}^\infty (z+KR/2)^{d-1}\bar{\alpha}(z)=O(K^d),$$
and $|d|=O(K^{2+d} |\D_n|m_n^{2d})$.

To bound $|a|$ we can combine the previous reasoning used to deal with
$|d|$ and the method used to bound $|a|$ in case
\ref{casetime}.
\begin{align*}
|a| &\leq 8 M_{\epsilon}^{4} \sum_{k=2}^{K} \sum_{l'=1}^{k-1}  \sum_{\substack{k'=1\\ l'-k'\le \lkn}}^{l'}\sum_{(\l,\j),(\l',\j')\in\mathcal{J}_{K,n}}  \Big( \tilde{\alpha}(k-1-l'; d(\{\l,\j\}, \{\l', \j'\})-(k'+1-l')R/2)\\
 &~~~~~~~~~~~~~~~~~~~~+I[k'=l']\sum_{i=1}^{k'-1}\tilde{\alpha}(k-1-i; d(\{\l,\j\}, \{\l', \j'\})-(k'+1-i)R/2)  \Big)\\
&\leq 8 M_{\epsilon}^{4}\sum_{(\l,\j),(\l',\j')\in\mathcal{J}_{K,n}}\left(l_K \sum_{k=2}^{K} \sum_{t=0}^{k-2}\tilde{\alpha}(t; d(\{\l,\j\}, \{\l', \j'\})-KR/2)\right.\\
&~~~~~~~~~~~~~~~~~~~~+\left.\sum_{k=2}^{K} \sum_{l'=1}^{k-1}\sum_{t=k-l'}^{k-2}\tilde{\alpha}(t; d(\{\l,\j\}, \{\l', \j'\})-KR/2) \right)\\
&\leq 8 M_{\epsilon}^{4}\sum_{(\l,\j),(\l',\j')\in\mathcal{J}_{K,n}}\left(l_KK \sum_{t=0}^{K}\tilde{\alpha}(t; d(\{\l,\j\}, \{\l', \j'\})-KR/2)\right.\\
&+\left. \sum_{t=0}^{K-2}\tilde{\alpha}(t; d(\{\l,\j\}, \{\l', \j'\})-KR/2)\#\{(k, l'): 1\leq l'<k\leq K, k-l'\leq t\leq k-2\} \right)\\
&\leq 8 M_{\epsilon}^{4}\sum_{(\l,\j),(\l',\j')\in\mathcal{J}_{K,n}}\left(l_KK \sum_{t=0}^{K}\tilde{\alpha}(t; d(\{\l,\j\}, \{\l', \j'\})-KR/2)\right.\\
&+\left. K\sum_{t=0}^{K-2} t\tilde{\alpha}(t; d(\{\l,\j\}, \{\l', \j'\})-KR/2) \right).
\end{align*}
Using the same reasoning as in the bound of $|d|$ and condition
\ref{cond:space-time conditional independence} we get
$$\sum_{(\l,\j),(\l',\j')\in\mathcal{J}_{K,n}}\tilde{\alpha}(t; d(\{\l,\j\}, \{\l', \j'\})-KR/2) = O(|\D_n|m_n^{2d} t^\psi K^{d}).$$
Hence
$$l_K K \sum_{t=0}^K\sum_{(\l,\j),(\l',\j')\in\mathcal{J}_{K,n}}\tilde{\alpha}(t; d(\{\l,\j\}, \{\l', \j'\})-KR/2)= O(l_K |\D_n| m_n^{2d}  K^{d+1})$$ 
and
$$\sum_{(\l,\j),(\l',\j')\in\mathcal{J}_{K,n}}\tilde{\alpha}(t; d(\{\l,\j\}, \{\l', \j'\})-KR/2)\underset{t\rightarrow\infty}{\longrightarrow} 0.$$
Hence, by Cesàro summation,
$$\frac{1}{K}
\sum_{t=0}^{K-2}\sum_{(\l,\j),(\l',\j')\in\mathcal{J}_{K,n}}t\tilde
\alpha(t; d(\{\l,\j\}, \{\l', \j'\})-KR/2) \underset{K\rightarrow\infty}{\longrightarrow} 0.$$
Thus $|a|=O(l_K |\D_n| m_n^{2d} K^{d+1})+o(K^2)$. In conclusion, $a,b,c,d$ are all $O(K^{2+d} |\D_n| m_n^{2d})+o(K^2)$.\\[\bsl]

In summary, $V_{K,n}^4 \EE_{\H_{0}}[ |A_1|^2 ]$ is
$O(\mkn^{2d}|\D_n|)$ in case \ref{casespace}, $o(K^2)$ in case
\ref{casetime}, and  $O(K^{2+d} |\D_n| m_n^{2d})+o(K^2)$ in case \ref{caseboth}.
By \eqref{eq:approxvar} and assumption~\ref{assumption variance
  liminf1}, \ref{assumption variance liminf2}, or \ref{assumption variance liminf3},
$|\D_n|/V_{K,n}^{2}$, $K/V_{K,n}^2$ or
$|\D_n|K/\sigma_{K,n}^2$
are $O(1)$ as $n$,  $K$, or both $n$ and $K$  tends to
infinity. Hence  $\EE_{\H_{0}} [|A_1|^2 ]$ is $O(\mkn^{2d}|\D_n|^{-1})$ in case
  \ref{casespace} and $o(1)$ in case \ref{casetime}. In case
  \ref{caseboth} we use $m_n=|\D_n|^\tau$ and $K=O(|\D_n|^C)$ to
  obtain that $\EE_{\H_{0}} [|A_1|^2 ]$ is $O(K^d  m_n^{2d}
  |\D_n|^{-1})+o(1)=O(|\D_n|^{d(C+2 \tau) -1})+o(1)$. Thus by
  \eqref{eq:mkn} and \eqref{eq:lkn} and since $C<1/d-2 \tau$, $\EE_{\H_{0}}[ A_1]$
converges to zero in all cases \ref{casespace}-\ref{caseboth}.

\subsubsection{Convergence of $\EE_{\H_0}[A_2]$}
According to Taylor's formula, there exists a constant $c$ so that
\[ \left | 1- it \bar{S}_{K,n}(k,\l) - \exp\left\{-it \bar{S}_{K,n}(k,\l)\right\} \right | \le c t^2 \bar S_{K,n}^2(k,\l).\]
We thus get
\begin{align*}
\EE_{\H_{0}} [\left| A_2 \right| ] &\leq \frac{ct^2}{V_{K,n}} \sum_{k=1}^K  \sum_{\l \in \D_n} \EE_{\H_{0}}\left[ |E_{k,n}(\l)| \bar S_{K,n}^2(k,\l) \right] \\
 &= \frac{ct^2}{V_{K,n}^{3}} \sum_{k=1}^K \sum_{\l \in \D_n}  \sum_{(l,\j), (l',\j') \in \mathcal{I}_{K,n}(k,\l)} \EE_{\H_{0}} \left[ |E_{k,n}(\l)| E_{l,n}(\j) E_{l',n}(\j')  \right].
\end{align*}
By conditional centering and since $l,l' \geq k$, $\EE_{\H_{0}} [|E_{k,n}(\l)| E_{l,n}(\j) E_{l',n}(\j') ]$ is zero
unless $l=l'$. By H{\" o}lder's inequality
\[
  \big| \EE_{\H_{0}} \left[|E_{k,n}(\l)| E_{l,n}(\j) E_{l,n}(\j') \right] \big| \leq \left( \EE_{\H_{0}} \left[|E_{k,n}(\l)|^3 \right] \EE_{\H_{0}} \left[| E_{l,n}(\j)|^{3} \right] \EE_{\H_{0}} \left[|E_{l,n}(\j')|^{3} \right] \right)^{\frac{1}{3}},
\]
and hence in case \ref{casespace} from \ref{cond:boundedmoment1}, we obtain
\begin{align*}
 & \EE_{\H_{0}} [\left| A_2 \right| ] \leq \frac{ct^2}{V_{K,n}^{3}}
                                  \sum_{k=1}^K \sum_{\l \in \D_n}
                                  \sum_{l=k}^K\sum_{\substack{\j,\j'
                                  \D_n\\ \dd(\l,\j)\le
  \mkn,\dd(\l,\j') \le \mkn}} \EE_{\H_{0}} \left[ |E_{k,n}(\l)|
  E_{l,n}(\j) E_{l,n}(\j')  \right]\\
     \le &     \frac{ct^2}{V_{K,n}^{3}}
                                  \sum_{k=1}^K \sum_{\l \in \D_n}
                                  \sum_{l=k}^K\sum_{\substack{\j,\j'\D_n\\ \dd(\l,\j)\le  \mkn,\dd(\l,\j') \le \mkn}} M_{\epsilon}^{3}
 \le   \frac{ct^2}{V_{K,n}^{3}}K^2|\D_n| (2\lceil \mkn \rceil + 1)^{2d}  M_{\epsilon}^{3},
\end{align*}
since the cardinality of $\{\j\in \Z^{d}: \dd(\l,\j) \leq \mkn\}$ is at most $(2\lceil \mkn \rceil + 1)^{d}$.\\
In the case \ref{casetime}, \ref{cond:boundedmoment2} implies that
\begin{align*}
\EE_{\H_{0}} [\left| A_2 \right| ] \leq \frac{ct^2}{V_{K,n}^{3}}  \sum_{k=1}^K \sum_{\l,\j,\j' \in \D_n} \sum_{\substack{k\leq l \leq K\\ l-k\leq \lkn} } M_{\epsilon}^3 \le  \frac{ct^2}{V_{K,n}^{3}} K |\D_n|^3 \lkn M_{\epsilon}^3.
\end{align*}
For case~\ref{caseboth} we obtain from \ref{cond:boundedmoment3} that
\begin{align*}
  \EE_{\H_{0}} [\left| A_2 \right| ] \leq \frac{ct^2}{V_{K,n}^{3}}
   \sum_{k=1}^K \sum_{\l \in \D_n}
   \sum_{\substack{k\leq l \leq K\\ l-k\leq \lkn} }
                                  \sum_{\substack{\j,\j'
                                  \D_n\\ \dd(\l,\j)\le
  \mkn,\dd(\l,\j') \le \mkn}} M_{\epsilon}^{3}
     \le     \frac{ct^2}{V_{K,n}^{3}}K|\D_n| (2\lceil \mkn \rceil + 1)^{2d} l_K M_{\epsilon}^{3}.
\end{align*}
Since $|\D_n|^{3/2}/V_{K,n}^3$, $K^{3/2}/V_{K,n}^3$, and
$|\D_n|^{3/2}K^{3/2}/V_{K,n}^3$
are $O(1)$ in respectively
cases~\ref{casespace}, \ref{casetime}, and \ref{caseboth} and using \eqref{eq:mkn} and \eqref{eq:lkn}, $\EE_{\H_{0}}[\left | A_2 \right | ]$ converges to zero in each of the cases~\ref{casespace}-\ref{caseboth}.

\subsubsection{Convergence of $\EE_{\H_0}[A_3]$}

Note that
\begin{align*}
  \big| \EE_{\H_{0}} [A_3] \big| \le & \frac{1}{V_{K,n}} \sum_{\l \in \D_n} \sum_{k=1}^K \left| \Cov_{\H_{0}}\left[E_{k,n}(\l), \exp\left\{it \left(\bar{S}_{K,n}- \bar{S}_{K,n}(k,\l) \right) \right\} \right] \right|.
\end{align*}
We can partition the difference $\bar{S}_{K,n}- \bar{S}_{K,n}(k,\l)$
as $B_1 + B_2$, where
\[
  B_1=\frac{1}{V_{K,n}}\sum_{l=1}^{k-1}  \sum_{\j \in \D_n}
E_{l,n}(\j)
\]
and

\[
  B_2= \begin{cases} \frac{1}{V_{K,n}} \sum_{l=k}^{K}
\sum_{\substack{\j \in \D_n: \dd(\l,\j) >  \mkn}} E_{l,n}(\j) & \text{in case \ref{casespace}} \\
 \frac{1}{V_{K,n}} \sum_{l=k: l-k>\lkn}^{K}
 \sum_{\j \in \D_n} E_{l,n}(\j) & \text{in case \ref{casetime}} \\
 \frac{1}{V_{K,n}} \left [\sum_{l=k: l-k \le \lkn}^{K}
 \sum_{\substack{\j \in \D_n: \dd(\l,\j) >  \mkn}} E_{l,n}(\j) + \sum_{l=k: l-k>\lkn}^{K}
       \sum_{\j \in \D_n} E_{l,n}(\j) \right ]
& \text{in case \ref{caseboth}.}
\end{cases}
\]

Since $B_1$ is $\H_{k-1}$ measurable and using conditional centering,
\begin{align*}
  |\Cov_{\H_{0}} & \left[ E_{k,n}(\l),  \exp\left\{it \left(\bar{S}_{K,n}-
        \bar{S}_{K,n}(k,\l) \right) \right\} \right]| = |\Cov_{\H_{0}}\left[
    E_{k,n}(\l),  \e^{it (B_1+B_2)}  \right] |\\
    &=|\EE_{\H_{0}}\left[ \e^{itB_1} \EE_{\H_{k-1}}\left[ E_{k,n}(\l)  \e^{itB_2}  \right]\right]| \le \EE_{\H_{0}} \left | \EE_{\H_{k-1}}\left[ E_{k,n}(\l)  \e^{itB_2}  \right]\right |.
\end{align*}

In case \ref{casespace}, $B_2$ is a function of $\{\Ekn{l,n}(\j): l=k, \ldots, K, \j\in\D_n, \dd(\l,\j) > \mkn\}$ and
\begin{align*}
    \EE_{\H_{k-1}}\left[ E_{k,n}(\l)  \e^{itB_2}  \right]= \EE_{\H_{k-1}}\left[ E_{k,n}(\l) \EE_{\H_{k}}\left[ \e^{itB_2} \right] \right]
    =  \Cov_{\H_{k-1}}\left[ E_{k,n}(\l),  \EE_{\H_{k}} \left[ \e^{it B_2} \right] \right]  .
\end{align*}
The localization assumption \ref{cond:spacelocalized} implies that $\EE_{\H_{k}}\left[ \e^{itB_2} \right]$ is $\H_{0,k,\cup_{\j \in\D_n: \dd(\l,\j) > \mkn} C(\j)}$ measurable  and hence by \eqref{eq:covalphabound3} and \ref{cond:spatial conditional independence},
\begin{align*}
&\EE_{\H_{0}} \left[ \left| \Cov_{\H_{k-1}}[ E_{k,n}(\l), \EE_{\H_k}[ \e^{it B_2}] ] \right |\right]    \le  8 M_{\epsilon} \EE_{\H_{0}}\left[ \alpha_{\H_{k-1}}\left(\H_{0,k, C(\l)}, \H_{0,k,\cup_{\j\in\D_n: \dd(\l,\j) > \mkn} C(\j)} \right)  \right]^{\frac{p-2}{p}}\\
    \le & 8 M_{\epsilon} \EE_{\H_{0}}\left[  \alpha_{K,R^2, \infty}\left(\dd(C(\l), \cup_{\j\in\D_n:
      \dd(\l,\j) > \mkn} C(\j)) \right)
      \right]^{\frac{p-2}{p}}
    \leq 8 M_{\epsilon} \bar{\alpha}(\mkn - 2R) \leq \text{ constant } \mkn^{\eta}.
\end{align*}

In case~\ref{casetime}, $B_2$ is a function of $\{\Ekn{l,n}(\j): l=k+ \lfloor \lkn \rfloor +1, \ldots, K, \j\in\D_n\}$ and
\begin{align*}
  \EE_{\H_{k-1}}\left[ E_{k,n} \e^{it B_2} \right] = \EE_{\H_{k-1}} \left[ E_{k,n}(\l) \EE_{\H_{k + \lfloor \lkn \rfloor }} \left[\e^{it B_2}  \right] \right]
  = \Cov_{\H_{k-1}} \left[ E_{k,n}(\l), \EE_{\H_{k + \lfloor \lkn \rfloor }}\left[\e^{it B_2}  \right]  \right] .
\end{align*}
By the localization assumption \ref{cond:timelocalized}, $\EE_{\H_{k + \lfloor \lkn \rfloor }}[\e^{itB_2} ]$ is $\H_{(k+\lkn-r)^{+}, k+\lfloor \lkn \rfloor, \D_n}$ measurable and using same type of arguments as for \ref{casespace},
\begin{align*}
  &\EE_{\H_{0}} \left[ \left | \EE_{\H_{k-1}}\left[ E_{k,n}(\l)  \e^{itB_2}  \right]\right | \right]  \le \EE_{\H_{0}} \left[ \left |\Cov_{\H_{k-1}} \left[ E_{k,n}(\l), \EE_{\H_{k + \lfloor \lkn \rfloor }}\left[\e^{it B_2}  \right]  \right]\right | \right]\\
    \le & 8 M_{\epsilon} \EE_{\H_0}\left[ \alpha_{\H_{k-1}}(\H_{0,k,\D_n}, \H_{(k + \lkn - r)^{+}, k + \lkn, \D_n}) \right]^{\frac{p-2}{p}}
    \le 8 M_{\epsilon} \EE_{\H_0}\left[ \rwrev{\alpha_r}(\lkn)  \right] ^{\frac{p-2}{p}} \\
    \le & 8 M_{\epsilon} \tilde{\alpha}(\lkn) \leq \text{ constant } \lkn^{\psi}.
\end{align*}

Considering case \ref{caseboth}, $B_2=B_{21}+B_{22}$ with
\[ B_{21} = \frac{1}{V_{K,n}} \sum_{l=k: l-k \le \lkn}^{K}
 \sum_{\substack{\j \in \D_n: \dd(\l,\j) >  \mkn}} E_{l,n}(\j)    \quad \text{and} \quad B_{22}=\frac{1}{V_{K,n}} \sum_{l=k: l-k>\lkn}^{K}
 \sum_{\j \in \D_n} E_{l,n}(\j) .
\]
Thus
\[ \left |\EE_{\H_{k-1}}\left[ E_{k,n}(\l)  \e^{itB_{2}}\right] \right | \le  \left |\EE_{\H_{k-1}}\left[ E_{k,n}(\l)\e^{itB_{22}}
\right] \right |+ \left |\EE_{\H_{k-1}}\left[ E_{k,n}(\l)
  (\e^{itB_{21}}-1) \e^{itB_{22}}\right] \right |.\]
Since $B_{22}$ coincides with $B_2$ in case \ref{casetime} we can
proceed as in case \ref{casetime} for the term $\EE_{\H_{k-1}} \!\!\left[ E_{k,n}(\l)\e^{itB_{22}}\right]$ in
the above inequality and conclude $\EE_{\H_0}\left[\left|\EE_{\H_{k-1}}\left[
  E_{k,n}(\l)\e^{itB_{22}}\right] \right |\right]$ is $O(l_K^\psi)$. For the second term, by a first order Taylor
expansion and (2) in the supplementary
material,
\begin{align*} &\left |\EE_{\H_{k-1}}\left[ E_{k,n}(\l)
   (\e^{itB_{21}}-1) \e^{itB_{22}}\right] \right | \le \sqrt{2} |t|
 \EE_{\H_{k-1}} \left[| E_{k,n}(\l)B_{21} |\right] \\ \le
               &\frac{\sqrt{2}|t|}{V_{k,n}} \sum_{l=k: l-k \le
                 \lkn}^{K}\sum_{\j \in \D_n} \EE_{\H_{k-1}}
                 \left[|E_{k,n}(\l) E_{l,n}(\j)|\right] = O(l_K |\D_n| M_{\epsilon}^2/V_{k,n}).
\end{align*}
We have $|\D_n|/V_{K,n}^{2}$, $K/V_{K,n}^2$ or
$|\D_n|K/\sigma_{K,n}^2$
are $O(1)$ as $n$,  $K$, or both $n$ and $K$  tends to
infinity. Using \eqref{eq:mkn} or \eqref{eq:lkn},  $\EE_{\H_{0}}[A_{3}]\to 0$ in cases \ref{casespace}-\ref{caseboth}.

\end{proof}



\section*{Acknowledgements}
We are grateful for the thoughtful and constructive comments from the
editors and reviewers that led to very substantial improvements of the
paper.  Rasmus Waagepetersen was supported by `urbanLab: Spatial data center for evidence-based city
  planning' (VIL57389), Villum Fonden, and by `NSECURE'
  (NNF23OC0084252), Novo Nordisk Foundation.

\section*{Supplementary material}
The supplementary material \rwrev{\citep{jalilian:poinas:xu:waagepetersen:23:suppl}} for this article provides supporting results and computations for the proof of \autoref{thm:clt} and for Section~\ref{sec:verifyspacetimecar}. It can be found at http://doi.org/10.1017/[TO BE SET].

\bibliographystyle{apalike} 
\bibliography{biblio_biscio}

\begin{thebibliography}{}

\bibitem[Besag, 1974]{besag1974spatial}
Besag, J. (1974).
\newblock Spatial interaction and the statistical analysis of lattice systems.
\newblock {\em J. R. Stat. Soc. Ser. B. Stat. Methodol.}, 36(2):192--225.

\bibitem[Biscio et~al., 2018]{biscio:poinas:waagepetersen:17}
Biscio, C. A.~N., Poinas, A., and Waagepetersen, R. (2018).
\newblock A note on gaps in proofs of central limit theorems.
\newblock {\em Statist. Probab. Lett.}, pages 7--10.

\bibitem[Biscio and Svane, 2022]{biscio:svane:22}
Biscio, C. A.~N. and Svane, A.~M. (2022).
\newblock {A functional central limit theorem for the empirical Ripley's
  K-function}.
\newblock {\em Electron. J. Stat.}, 16(1):3060--3098.

\bibitem[Biscio and Waagepetersen, 2019]{biscio:waagepetersen:18}
Biscio, C. A.~N. and Waagepetersen, R. (2019).
\newblock A general central limit theorem and a subsampling variance estimator
  for $\alpha$-mixing point processes.
\newblock {\em Scand. J. Stat.}, 46(4):1168--1190.

\bibitem[Bolthausen, 1982]{bolthausen:82}
Bolthausen, E. (1982).
\newblock On the central limit theorem for stationary mixing random fields.
\newblock {\em Ann. Probab.}, 10(4):1047--1050.

\bibitem[Bradley, 2005]{bradley2005basic}
Bradley, R.~C. (2005).
\newblock Basic properties of strong mixing conditions. a survey and some open
  questions.
\newblock {\em Probab. Surv.}, 2:107--144.

\bibitem[Brown, 1971]{brown1971martingale}
Brown, B.~M. (1971).
\newblock Martingale central limit theorems.
\newblock {\em Ann. Math. Stat.}, 42:59--66.

\bibitem[Cheysson and Lang, 2022]{cheysson2022spectral}
Cheysson, F. and Lang, G. (2022).
\newblock Spectral estimation of {H}awkes processes from count data.
\newblock {\em The Annals of Statistics}, 50(3):1722--1746.

\bibitem[Comets and Jan{\c z}ura, 1998]{comets:janzura:98}
Comets, F. and Jan{\c z}ura, M. (1998).
\newblock A central limit theorem for conditionally centred random fields with
  an application to {M}arkov fields.
\newblock {\em J. Appl. Probab.}, 35(3):608--621.

\bibitem[Condit et~al., 2019]{condit:etal:19}
Condit, R., Perez, R., Aguilar, S., Lao, S., Foster, R., and Hubbell, S.
  (2019).
\newblock Complete data from the barro colorado 50-ha plot: 423617 trees, 35
  years [dataset]., 2019 version.
\newblock Dryad, https://doi.org/10.15146/5xcp-0d46.

\bibitem[Dahlhaus, 2000]{dahlhaus2000graphical}
Dahlhaus, R. (2000).
\newblock Graphical interaction models for multivariate time series.
\newblock {\em Metrika}, 51(2):157--172.

\bibitem[Davydov, 1968]{davydov:68}
Davydov, Y.~A. (1968).
\newblock Convergence of distributions generated by stationary stochastic
  processes.
\newblock {\em Theory Probab. Appl.}, 13(4):691--696.

\bibitem[Doukhan, 1994]{doukhan:94}
Doukhan, P. (1994).
\newblock {\em Mixing, Properties and Examples}, volume~85 of {\em Lecture
  Notes in Statistics}.
\newblock Springer-Verlag, New York.

\bibitem[Doukhan and Lang, 2016]{doukhan:lang:16}
Doukhan, P. and Lang, G. (2016).
\newblock Weak dependence of point processes and application to second-order
  statistics.
\newblock {\em Statistics}, 50(6):1221--1235.

\bibitem[Guyon, 1995]{guyon:95}
Guyon, X. (1995).
\newblock {\em Random Fields on a Network}.
\newblock Probability and its Applications. Springer, New York.

\bibitem[Hall et~al., 2014]{hall:heyde:14}
Hall, P., Heyde, C., Birnbaum, Z., and Lukacs, E. (2014).
\newblock {\em Martingale Limit Theory and Its Application}.
\newblock Communication and Behavior. Elsevier Science.

\bibitem[Heinrich, 2013]{Heinrich2013}
Heinrich, L. (2013).
\newblock Asymptotic methods in statistics of random point processes.
\newblock In Spodarev, E., editor, {\em Stochastic Geometry, Spatial Statistics
  and Random Fields: Asymptotic Methods}, pages 115--150. Springer Berlin
  Heidelberg, Berlin, Heidelberg.

\bibitem[Heinrich and Pawlas, 2013]{heinrich2013absolute}
Heinrich, L. and Pawlas, Z. (2013).
\newblock Absolute regularity and brillinger-mixing of stationary point
  processes.
\newblock {\em Lithuanian Mathematical Journal}, 53(3):293--310.

\bibitem[Jalilian et~al.,
  2024]{jalilian2024compositelikelihoodinferencespacetime}
Jalilian, A., Cuevas-Pacheco, F., Xu, G., and Waagepetersen, R. (2024).
\newblock Composite likelihood inference for space-time point processes.
\newblock arXiv 2402.12548.

\bibitem[Jalilian et~al., 2023]{jalilian:poinas:xu:waagepetersen:23:suppl}
Jalilian, A., Poinas, A., Xu, G., and Waagepetersen, R. (2023).
\newblock Supplementary material for `a central limit theorem for a sequence of
  conditionally centered random fields'.
\newblock Supplementary material.

\bibitem[Jensen and K{\"u}nsch, 1994]{jensen1994}
Jensen, J.~L. and K{\"u}nsch, H.~R. (1994).
\newblock On asymptotic normality of pseudo likelihood estimates for pairwise
  interaction processes.
\newblock {\em Ann. Inst. Statist. Math.}, 46:475--486.

\bibitem[Kar{\'a}csony, 2006]{karacsony:06}
Kar{\'a}csony, Z. (2006).
\newblock A central limit theorem for mixing random fields.
\newblock {\em Miskolc Math. Notes}, 7(2):147--160.

\bibitem[Leonenko, 1975]{leonenko:75}
Leonenko, N.~N. (1975).
\newblock Central limit theorem for $m$-dependent random fields in schemes
  related to series schemes.
\newblock {\em Ukrainian Math. J.}, 27:556--559.

\bibitem[Liu et~al., 2024]{liu2024two}
Liu, W., Xu, G., Fan, J., and Zhu, X. (2024).
\newblock Two-way homogeneity pursuit for quantile network vector
  autoregression.
\newblock arXiv 2404.18732.

\bibitem[Poinas et~al., 2019]{poinas:delyon:lavancier:19}
Poinas, A., Delyon, B., and Lavancier, F. (2019).
\newblock {Mixing properties and central limit theorem for associated point
  processes}.
\newblock {\em Bernoulli}, 25(3):1724--1754.

\bibitem[Prakasa~Rao, 2009]{rao:09}
Prakasa~Rao, B. (2009).
\newblock Conditional independence, conditional mixing and conditional
  association.
\newblock {\em Ann. Inst. Statist. Math.}, 61:441--460.

\bibitem[Rio, 1993]{rio:93}
Rio, E. (1993).
\newblock Covariance inequalities for strongly mixing processes.
\newblock {\em Ann. Inst. Henri Poincar{\'e} Probab. Stat.}, 29(4):587--597.

\bibitem[Waagepetersen and Guan, 2009]{waagepetersen:guan:09}
Waagepetersen, R. and Guan, Y. (2009).
\newblock Two-step estimation for inhomogeneous spatial point processes.
\newblock {\em J. R. Stat. Soc. Ser. B. Stat. Methodol.}, 71(3):685--702.

\bibitem[Yuan and Lei, 2013]{yuan2013some}
Yuan, D. and Lei, L. (2013).
\newblock Some conditional results for conditionally strong mixing sequences of
  random variables.
\newblock {\em Sci. China Math.}, 56(4):845--859.

\bibitem[Zhu et~al., 2017]{zhu2017network}
Zhu, X., Pan, R., Li, G., Liu, Y., and Wang, H. (2017).
\newblock Network vector autoregression.
\newblock {\em Ann. Statist.}, 45:1096--1123.

\bibitem[Zhu et~al., 2023]{zhu2023simultaneous}
Zhu, X., Xu, G., and Fan, J. (2023).
\newblock Simultaneous estimation and group identification for network vector
  autoregressive model with heterogeneous nodes.
\newblock {\em Journal of Econometrics}, page 105564.

\end{thebibliography}


\begin{thebibliography}{}

\bibitem[Biscio and Waagepetersen, 2019]{biscio:waagepetersen:18}
Biscio, C. A.~N. and Waagepetersen, R. (2019).
\newblock A general central limit theorem and a subsampling variance estimator
  for $\alpha$-mixing point processes.
\newblock {\em Scand. J. Stat.}, 46(4):1168--1190.

\bibitem[Hamilton, 1994]{hamilton1994time}
Hamilton, J.~D. (1994).
\newblock {\em Time Series Analysis}.
\newblock Princeton University Press.

\bibitem[Jalilian et~al., 2023]{jalilian:poinas:xu:waagepetersen:23}
Jalilian, A., Poinas, A., Xu, G., and Waagepetersen, R. (2023).
\newblock A central limit theorem for a sequence of conditionally centered
  random fields.
\newblock Submitted.

\bibitem[Rencher and Schaalje, 2008]{rencher2008linear}
Rencher, A.~C. and Schaalje, G.~B. (2008).
\newblock {\em Linear Models in Statistics}.
\newblock John Wiley \& Sons.

\bibitem[Yamamoto, 1981]{yamamoto1981predictions}
Yamamoto, T. (1981).
\newblock Predictions of multivariate autoregressive-moving average models.
\newblock {\em Biometrika}, 68(2):485--492.

\end{thebibliography}

\end{document}


\maketitle

\begin{abstract}
This document contains supplementary material for the paper `A central
limit theorem for a sequence of conditionally centered random fields'.

\textbf{Keywords:} Central limit theorem;  conditional $\alpha$-mixing; conditional centering;
estimating function; random field; space-time; spatio-temporal; vector autoregression\\
\textbf{MSC2020 subject classifications:} 60F05; 60G55; 60G60; 62M30
\end{abstract}


The following sections contain results on moments, mixing
coefficients, and covariances of relevance for the paper
\cite{jalilian:poinas:xu:waagepetersen:23}. Cross-references starting with
`M-' refer to that paper to which we also refer for notational
definitions. In Section~\ref{sec:boundsmoments} and
  Section~\ref{sec:boundscov} we consider the one-dimensional case $q=1$.

\section{Bounds on moments}\label{sec:boundsmoments}

If assumption M-($\mathcal{A}$4)
holds, applications of
Lyapunov's and H\"{o}lder's inequalities imply that
\begin{align}
  \EE_{\H_{0}} \left[ \left|\Ekn{k,n}(\l)\right|^{p} \right]  \leq \left( \EE_{\H_{0}} \left[ |\Ekn{k,n}(\l)|^{6+\epsilon} \right] \right)^{\frac{p}{6+\epsilon}} \leq \left( d_{K,n} \right)^{\frac{p}{6+\epsilon}}  \leq M_{\epsilon}^{p} \label{eq:p new moment bound1}
\end{align}
if $1\leq p \leq 6+\epsilon$,
\begin{align}
  \EE_{\H_{0}} & \left[ \left|\Ekn{k,n}(\l) \Ekn{l,n}(\j) \right|^{p} \right]  \leq \left( \EE_{\H_{0}}\left[ |\Ekn{k,n}(\l)|^{2p} \right]  \EE_{\H_{0}}\left[ |\Ekn{l,n}(\j)|^{2p} \right] \right)^{\frac{1}{2}} \nonumber\\
  &\leq   \left( M_{\epsilon}^{4p} \right)^{\frac{1}{2}}  \leq M_{\epsilon}^{2p} \label{eq:p new moment bound2}
\end{align}
if $1\leq p\leq 3 + \epsilon/2$, and
\begin{align}
  \EE_{\H_{0}} & \left[ \left|\Ekn{k,n}(\l) \Ekn{l,n}(\j) \Ekn{l',n}(\j') \right|^{p} \right]  \nonumber\\
  &\leq \left( \EE_{\H_{0}}\left[ |\Ekn{k,n}(\l)|^{3p} \right] \right)^{\frac{1}{3}} \left( \EE_{\H_{0}}\left[ |\Ekn{l,n}(\j)\Ekn{l',n}(\j')|^{\frac{3p}{2}} \right] \right)^{\frac{2}{3}} \nonumber\\
  &\leq \left( \EE_{\H_{0}}\left[ |\Ekn{k,n}(\l)|^{3p} \right]  \EE_{\H_{0}}\left[ |\Ekn{l,n}(\j)|^{3p} \right] \EE_{\H_{0}}\left[|\Ekn{l',n}(\j')|^{3p} \right] \right)^{\frac{1}{3}} \nonumber\\
  &\leq   \left( M_{\epsilon}^{9p} \right)^{\frac{1}{3}}  \leq M_{\epsilon}^{3p} \label{eq:p new moment bound3}
\end{align}
if $1\leq p \leq 2 + \epsilon/3$.
Further, by  Jensen's inequality, for any $1\leq k'\leq k\leq l\leq
K$,
\begin{align}
  \EE_{\H_{0}}\left[ \big| \EE_{\H_{k'-1}}\left[ \Ekn{k,n}(\l) \Ekn{l,n}(\j) \right] \big|^{p}  \right] &\leq \EE_{\H_{0}}\left[ \EE_{\H_{k'-1}}\left[ | \Ekn{k,n}(\l) \Ekn{l,n}(\j) |^{p}  \right]  \right] \leq M_{\epsilon}^{2p}. \label{eq:p new moment bound4}
\end{align}

\section{Monotonicity of the $\alpha$-mixing coefficients}\label{sec:mixingcoefficient}

If $m \leq m_2$, then
\[
  \{A,B\subset \D: \dd(A, B) \geq m_2 \} \subseteq   \{A,B\subset \D: \dd(A, B) \geq m \}
\]
Also , if $m' \leq m'_2$, then
\[
  \{l,k\in \Z_{+}: 0\leq l, k\leq K, l+m'_2 \leq k \} \subseteq \{l,k\in \Z_{+}: 0\leq l, k\leq K, l+m' \leq k \}
\]
hence
\[ \alpha_{K,c_1,c_2}(m_2) \leq \alpha_{K,c_1,c_2}(m),~~\alpha_{\rwrev{r}}(m'_2) \leq \alpha_{\rwrev{r}}(m') \]
and
\[ \apMOD{\alpha_{r; K,c_1,c_2}(m_2; m'_2) \leq \alpha_{r; K,c_1,c_2}(m; m')}. \]
Using the same arguments, for $c_{1} \leq c_{1}'$, $c_{2} \leq c_{2}'$ and $r \leq r'$ we get
\begin{align*}
\alpha_{K,c_1,c_2}(m) \leq \alpha_{K,c_{1}',c_2}(m),~~ \alpha_{K,c_1,c_2}(m) \leq \alpha_{K,c_{1},c_{2}'}(m),~~\alpha_{\rwrev{r}}(m) \leq \alpha_{\rwrev{r'}}(m)
\end{align*}
and
$$\apMOD{\alpha_{r; K,c_1,c_2}(m; m') \leq \alpha_{r'; K,c_{1}',c'_2}(m; m')}.$$

\section{Bounds on covariances}\label{sec:boundscov}

In this section we provide bounds for covariances of the form
\begin{align} \label{eq:covarianceofA1terms}
\Cov_{\H_{0}} & \left[ E_{k,n}(\l) E_{l,n}(\j), E_{k',n}(\l') E_{l',n}(\j') \right]
\end{align}
for $1\leq k,k'\leq K$ and $k\leq l\leq K, k'\leq l'\leq K$ in the
one-dimensional case $q=1$.

The covariance \eqref{eq:covarianceofA1terms}
 can be rewritten as
\begin{align*}
 &\EE_{\H_{0}} \left[ \Cov_{\H_{k^{*}-1}}\left[ E_{k,n}(\l) E_{l,n}(\j), E_{k',n}(\l') E_{l',n}(\j')  \right] \right] \\
   &+  \Cov_{\H_{0}}\left[ \EE_{\H_{k^{*}-1}} [ E_{k,n}(\l) E_{l,n}(\j) ], \EE_{\H_{k^*-1}}[ E_{k',n}(\l') E_{l',n}(\j') ]  \right],
\end{align*}
where $k^{*} = \max\{k,l,k',l'\}$.
If precisely one index coincides  with $k^*$, then by conditional centering, the covariance~\eqref{eq:covarianceofA1terms} is
zero. Thus, the covariance~\eqref{eq:covarianceofA1terms} can only be
non-zero in the cases
a) $k^*=k=l>l'\geq k'$ or $k^*=k'=l'>l \geq k$,
b) $k^*=l=l'>k,k'$,
c) $k^*=k=l=l'>k'$ or $k^*=l'=k'=l>k$, and
d) $k^*=k=l=k'=l'$.
Case a) is considered in Section~M-5.2.1 
of \cite{jalilian:poinas:xu:waagepetersen:23}. Cases b)-d) are considered for M-(i) and M-(iii) in the following subsections. We consider the spatial case M-(i) 
and the spacetime case M-(iii) 
at the same time since the bounds we require are identical for these cases and can be retrieved in a similar manner. For the time case M-(ii) 
we only require a constant bound. For b) and c) this can be recovered from the spacetime case by bounding the $\alpha$-mixing coefficient by its maximum value of $1/4$. For d) this can be done directly by bounding the covariance~\eqref{eq:covarianceofA1terms} with Hölder's inequality. The cases and associated bounds are summarized in M-Table~1 
in \cite{jalilian:poinas:xu:waagepetersen:23}.

\subsection{Case b)}

Consider $k^*=l=l'> k_{*}=\max\{k,k'\}$. Then
\[
  \EE_{\H_{l-1}}[ E_{k,n}(\l) E_{l,n}(\j) ] = \EE_{\H_{l-1}}[ E_{k',n}(\l') E_{l,n}(\j') ] = 0
\]
and hence the covariance~\eqref{eq:covarianceofA1terms} becomes
$\EE_{\H_{0}} \left[ \Cov_{\H_{l-1}}  \left[ E_{k,n}(\l) E_{l,n}(\j),
    E_{k',n}(\l') E_{l,n}(\j')  \right]  \right]$.

If either condition~M-($\mathcal{A}$2-1) 
or~M-($\mathcal{A}$2-3) 
holds, $E_{k,n}(\l)E_{l,n}(\j)$ and $E_{k,n}(\l')E_{l,n}(\j')$ are
respectively $\H_{0, l, C(\l)\cup C(\j)}$ and $\H_{0, l, C(\l')\cup
  C(\j')}$ measurable and  inequalities~\eqref{eq:p new moment bound2}
and~M-(8) 
imply that
\begin{align*}
& \left| \EE_{\H_{0}}\left[ \Cov_{\H_{l-1}}  \left[ E_{k,n}(\l) E_{l,n}(\j), E_{k',n}(\l') E_{l,n}(\j')\right] \right] \right|\\
 \leq & 8 M_{\epsilon}^{4} \EE_{\H_{0}}\left[\alpha_{\H_{l-1}}(\H_{0,l, C(\l)\cup C(\j)}, \H_{0,l, C(\l')\cup C(\j')})\right]^{\frac{p-2}{p}}.
\end{align*}
Thus, since
$$\alpha_{\H_{l-1}} \left(\H_{0, l, C(\l)\cup C(\j)}, \H_{0, l, C(\l')\cup C(\j')}\right)\leq \apMOD{\alpha_{l; K,2R^2, 2R^2}}\left(0; \dd(\{\l, \j\}, \{\l', \j'\}) - R \right)$$
from either condition~M-($\mathcal{A}$3-1) 
or condition~M-($\mathcal{A}$3-3) 
we obtain
\begin{align*}
   & \left| \EE_{\H_{0}}\left[\Cov_{\H_{l-1}}  \left[ E_{k,n}(\l) E_{l,n}(\j), E_{k',n}(\l') E_{l,n}(\j')\right] \right] \right|\\
  \leq & 8 M_{\epsilon}^{4}  \EE_{\H_{0}} \left[\apMOD{\alpha_{l; K,2R^2, 2R^2}}\left(0; \dd(\{\l, \j\}, \{\l', \j'\}) - R \right)\right]^{\frac{p-2}{p}}\\
  \leq & 8 M_{\epsilon}^{4} \bar{\alpha}(\dd(\{\l, \j\}, \{\l', \j'\}) - R).
\end{align*}

\subsection{Case c)}

Consider $k^*=k=l=l'>k_{*}=k'$. Then $\EE_{\H_{k-1}}[ E_{k,n}(\l)
E_{k',n}(\l') ]= 0$ and hence the
covariance~\eqref{eq:covarianceofA1terms} is given by $\EE_{\H_{0}}
\left[ \Cov_{\H_{k-1}}  \left[ E_{k,n}(\l) E_{k,n}(\j), E_{k',n}(\l')
    E_{k,n}(\j')  \right] \right]$. If either conditions~M-($\mathcal{A}$2-1) 
and M-($\mathcal{A}$3-1) 
or conditions~M-($\mathcal{A}$2-3) 
and M-($\mathcal{A}$3-3) 
hold, by inequalities~\eqref{eq:p new moment bound2} and M-(8),
\begin{align*}
 \big| \EE_{\H_{0}} & \left[ \Cov_{\H_{k-1}}  \left[ E_{k,n}(\l) E_{k,n}(\j), E_{k',n}(\l') E_{k,n}(\j')\right]  \right] \big| \\
 &\leq 8 M_{\epsilon}^{4} \EE_{\H_{0}} \left[ \alpha_{\H_{k-1}}(\H_{0,k,C(\l)\cup C(\j)}, \H_{0,k, C(\l')\cup C(\j')}) \right] ^{\frac{p-2}{p}} \\
  &\leq 8 M_{\epsilon}^{4} \EE_{\H_{0}} \left[  \apMOD{\alpha_{k; K,2R^2, 2R^2}}\left(0; \dd(\{\l,\j\},\{\l',\j'\})-R\right)\right] ^{\frac{p-2}{p}}\\
  &\leq 8 M_{\epsilon}^{4} \bar{\alpha}(\dd(\{\l,\j\},\{\l',\j'\}-R).
\end{align*}

\subsection{Case d)}
When $k^*=k_{*}=k=l=k'=l'$, using the same recursion idea as in case M-a), the covariance~\eqref{eq:covarianceofA1terms} can be rewritten as
\begin{align*}
  & \EE_{\H_{0}} \left[ E_{k,n}(\l) E_{k,n}(\j) \left(E_{k,n}(\l') E_{k,n}(\j') - \EE_{\H_{0}} [E_{k,n}(\l') E_{k,n}(\j')]\right) \right] \\
= & \sum_{i=1}^k \EE_{\H_{0}} \left[ E_{k,n}(\l) E_{k,n}(\j) \left(\EE_{\H_i} [E_{k,n}(\l') E_{k,n}(\j')] - \EE_{\H_{i-1}} [E_{k,n}(\l') E_{k,n}(\j')]\right) \right] \\
= & \sum_{i=1}^k \EE_{\H_{0}} \left[ \Cov_{\H_{i-1}}\left[ E_{k,n}(\l) E_{k,n}(\j), \EE_{\H_i} [E_{k,n}(\l') E_{k,n}(\j')]\right] \right].
\end{align*}
If either condition~M-($\mathcal{A}$2-1) 
or~M-($\mathcal{A}$2-3) 
holds, $E_{k,n}(\l) E_{k,n}(\j)$ and $\EE_{\H_i}\left[E_{k,n}(\l') E_{k,n}(\j')  \right]$ are respectively measurable with respect to $\H_{0,k,C(\l)\cup C(\j)}$ and $\H_{0,i,(k-i-1)(C(\l')\cup C(\j'))}$
Thus, since
\begin{align*}
  \alpha_{\H_{i-1}} & \left(\H_{0, k, C(\l)\cup C(\j)}, \H_{0, i, (k-i-1)(C(\l')\cup C(\j'))}\right) \\
   &\leq \apMOD{\alpha_{i+1; k, 2R^2, 2((k-i-1)R)^2}}\big(0; \dd(C(\l)\cup C(\j), (k-i-1)(C(\l')\cup C(\j'))) \big) \\
    &\leq \apMOD{\alpha_{i+1; k, 2R^2, \infty}}\left(0 ; \dd(\{\l, \j\}, \{\l', \j'\}) - (k-i)R/2\right),
\end{align*}
from inequalities~\eqref{eq:p new moment bound2}, M-(8)
and either condition~M-($\mathcal{A}$3-1) 
or condition~M-($\mathcal{A}$3-3) 
we obtain that \eqref{eq:covarianceofA1terms} is bounded by
\begin{align*}
   & \sum_{i=1}^{k} 8 M_{\epsilon}^{4} \EE_{\H_0}\left[ \apMOD{\alpha_{i+1; k,2R^2, \infty}}\left(0 ; \dd(\{\l, \j\}, \{\l', \j'\}) - (k-i)R/2 \right) \right] ^{\frac{p-2}{p}} \\
  \leq & 8 M_{\epsilon}^{4} \sum_{i=1}^{k} \bar \alpha(\dd(\{\l, \j\}, \{\l', \j'\}) - (k-i)R/2)\\
\leq & 8 M_{\epsilon}^{4} K \bar \alpha(\dd(\{\l, \j\}, \{\l', \j'\}) - KR/2).
\end{align*}

\section{Convergence of variance}\label{sec:convvariance}

In the case M-(ii), 
$m_{n}=\infty$ and hence $\mathcal{J}_{K,n}=\D_n\times \D_n$  and $\sigma^{2}_{K,n} = V^{2}_{K,n}$.
But in cases M-(i) and M-(iii), 
\begin{align*}
  \sigma^{2}_{K,n} &= V^{2}_{K,n} +  \sum_{k=1}^K \sum_{\substack{\l,\j\in \D_n \\ \dd(\l,\j) > m_{n}}} \EE_{\H_{0}}\left[ E_{k,n}(\l) E_{k,n}(\j) \right]
\end{align*}
and for any $k,n\in\N$, $\l,\l',\j,\j'\in\Z^{d}$ and $2<p\leq 2 + \epsilon/3$,  from  inequalities~\eqref{eq:p new moment bound1} and M-(8), 
\begin{align}
 &\big | \EE_{\H_{0}}  \left[ E_{k,n}(\l) E_{k,n} (\j) \right]  \big| = \big|  \EE_{\H_0} \left[\Cov_{\H_{k-1}}  \left[ E_{k,n}(\l), E_{k,n} (\j) \right]\right] \big| \nonumber\\
  \le &8 \EE_{\H_0}\left[ \left| E_{k,n}(\l) \right|^{p}
    \right]^{\frac{1}{p}}\EE_{\H_0}\left[ \left| E_{k,n}(\j)
    \right|^{p}  \right]^{\frac{1}{p}} \EE_{\H_0} \left[
    \alpha_{\H_{k-1}} \big( \H_{0,k,C(\l)}, \H_{0,k,C(\j)} \big)
    \right] ^{\frac{p-2}{p}}\nonumber \\\le &8M_\epsilon^2 \EE_{\H_0} \left[\alpha_{\H_{k-1}} \big( \H_{0,k,C(\l)}, \H_{0,k,C(\j)} \big)
    \right] ^{\frac{p-2}{p}} . \nonumber
\end{align}

Thus, in case (i), 
using M-($\mathcal{A}$3-1) 
\begin{align*}
  |\sigma^{2}_{K,n} - V^{2}_{K,n}| &\leq \sum_{k=1}^K \sum_{\substack{\l,\j\in \D_n \\ \dd(\l,\j) > m_{n}}} \left | \EE_{\H_{0}} \left [ \Cov_{\H_{k-1}}[ E_{k,n}(\l), E_{k,n}(\j) ] \right] \right |\\
  &\leq 8 M_{\epsilon}^{2} \sum_{k=1}^K \sum_{\substack{\l,\j\in \D_n \\ \dd(\l,\j) > m_{n}}}  \EE_{\H_{0}}\left[\alpha_{\H_{k-1}}(\H_{0, k, C(\l)}, \H_{0, k, C(\j)})  \right]^{\frac{p-2}{p}} \\
  &\leq 8 K M_{\epsilon}^{2} \sum_{\l \in \D_n} \sum_{\substack{\j \in \Z^{d} \\ \dd(\l, \j) > m_{n}}} \EE_{\H_{0}}\left[ \alpha_{K,R^2, R^2}(\dd(\l,\j) - R) \right]^{\frac{p-2}{p}} \\
  &\leq 8 K M_{\epsilon}^{2} \sum_{\l \in \D_n} \sum_{l=\lfloor m_{n} \rfloor}^{\infty}  \sum_{\substack{\j\in \Z^{d} \\ \dd(\l,\j) = l}}  \bar{\alpha}(l - R) \\
  &\leq 8 (3^{d}) K |\D_n| M_{\epsilon}^{2} \sum_{l=\lfloor m_{n} \rfloor}^{\infty} l^{d-1} \bar{\alpha}(l - R)
\end{align*}
because the cardinality of $\{\l'\in\Z^{d}: \dd(\l,\l')=l\}$ is less than or equal to $3^{d} l^{d - 1}$ for $l\in \N$ (see \cite{biscio:waagepetersen:18}, Lemma~F4). Thus
\[
 \left|1 - \frac{V_{K,n}^2}{\sigma_{K,n}^2} \right| \leq 8 (3^{d}) K \frac{|\D_n|}{\sigma^{2}_{K,n}} M_{\epsilon}^{2} \sum_{l=\lfloor m_{n} \rfloor}^{\infty} l^{d-1} \bar{\alpha}(l - R).
\]
As $n\to \infty$, assumption~M-($\mathcal{A}$5-1) 
implies that $|\D_n|/\sigma_{K,n}^2$ is bounded.
Condition M-($\mathcal{A}$3-1) 
implies that
\[
  \sum_{l=\lfloor m_{n} \rfloor}^{\infty} l^{d-1} \bar{\alpha}(l - R) \leq \text{ constant } \sum_{l=\lfloor m_{n} \rfloor}^{\infty} l^{\eta + d-1}
\]
and since $\eta + d < 0$, the tail series $\sum_{l=\lfloor m_{n} \rfloor}^{\infty} l^{\eta + d-1}$ converges to zero as $n\to\infty$. Therefore
\begin{equation}\label{eq:V and sigma convergence}
  \left|1 - \frac{V_{K,n}^2}{\sigma_{K,n}^2} \right| \xrightarrow[n\rightarrow \infty ]{a.s.} 0.
\end{equation}
The case M-(iii) is handled analogously but using
  assumptions~M-($\mathcal{A}$3-3) and M-($\mathcal{A}$5-3).

\section{Spatio-temporal autoregressive model}

In this section we provide details regarding
Section~M-4.2. 

\subsection{Representation as Vector Autoregression}\label{sec:VAR}

By rewriting~M-(6) 
as a Vector Autoregression (VAR) equation of order one (see
\cite{hamilton1994time}, p.~259) and  recursively applying the
resulting $\mathrm{VAR}(1)$ representation, it follows that for any $0\leq k < l$,
\begin{equation}\label{eq:recursive Xkn}
    X_{l} = C_{l-k} X_{k} + \sum_{j=0}^{l-k-1} C_{j} Q^{-1/2} Z_{l-j}
\end{equation}
where $C_{j} = H^{\T} F^{j} H$ and (see \cite{yamamoto1981predictions}, Lemma~1)
\[
  F=\begin{bmatrix}
  \beta B_{1} & \beta B_{2} & \cdots & \beta B_{r-1} & \beta B_{r} \\
  I_{|\D|} & 0 & \cdots & 0 & 0 \\
  0 & I_{|\D|} & \cdots & 0 & 0 \\
  \vdots & \vdots & \ddots & \vdots & \vdots \\
  0 & 0 & \cdots & I_{|\D|} & 0
  \end{bmatrix}_{r|\D|\times r|\D|},
  \quad
  H = \begin{bmatrix}
  I_{|\D|} \\
  0 \\
  0 \\
  \vdots \\
  0
  \end{bmatrix}_{r|\D|\times |\D|}.
\]
Thus, for any $0\leq k < l$, the conditional distribution of $X_{l}$ given $\H_{k}$ is multivariate normal with $\EE_{\H_{k}}[X_{l}] = C_{l-k} X_{k}$ and
\[
  \Var_{\H_{k}}[ X_{l}] = \sum_{j=0}^{l-k-1} C_{j} Q^{-1} C_{j}^{\T}.
\]
Moreover, for any $1\leq k \leq l$ and $m\in {\N}$,
\begin{align*}
  & \Cov_{\H_{k-1}}  \left[ X_{l}, X_{l+m} \right] =  
     \Cov_{\H_{k-1}}\left[ \EE_{\H_{l}}[X_{l}], \EE_{\H_{l}}[X_{l+m}] \right] \\
 =~&\Cov_{\H_{k-1}}\left[ X_{l}, C_{m}X_{l} \right] =\Cov_{\H_{k-1}}\left[ X_{l}, X_{l} \right] C_{m}^{\T}\\ = &\left(\sum_{j=1}^{l-k} C_{j} Q^{-1} C_{j}^{\T} \right) C_{m}^{\T}.
\end{align*}

If the eigenvalues of $F$ have modulus less than one, or equivalently the roots of the equation
\[
  \mathrm{det}\left(I_{|\D|} - \sum_{j=1}^{r} z^{j} \beta B_{j} \right) = 0
\]
are outside the complex unit circle, then $F^{m}\to0$  as $m\to\infty$(see \cite{hamilton1994time}, p.~259). We assume that as $m\to \infty$, $\| F^{m}\| = O(m^{\psi})$ for some $\psi<-1$, where for a matrix $A=[a_{i,j}]_{n_1\times n_2}$,
\[
  \| A\| = \left( \sum_{i=1}^{n_1}\sum_{j=1}^{n_2} a_{i,j}^{2} \right)^{\frac{1}{2}}
\]
is the Frobenius norm. Then, $\| C_{m}\|=O(m^{\psi})$ because
\[
  \| C_{m}\| = \|H^{\T} F^{m} H\| \leq \|H^{\T}\| \|F^{m}\| \| H\| \leq |\D|^2 \|F^{m}\|.
\]
Consequently for any $1\leq k \leq l$  and $\l,\j\in\D$,
\begin{equation}\label{eq:condcovboundautoreg}
  \left| \Cov_{\H_{k-1}} \left[ X_{l}(\l), X_{l+m}(\j) \right] \right| = O(m^{\psi}),
\end{equation}
as $m\to\infty$, since
\[
  \left\| \Cov_{\H_{k-1}} \left[ X_{l}, X_{l+m} \right] \right\| \leq \left\| \sum_{j=1}^{l-k} C_{j} Q^{-1} C_{j}^{\T}  \right\| \|C_{m}\|
\]
and
\[
  \left\| \sum_{j=1}^{l-k} C_{j} Q^{-1} C_{j}^{\T}  \right\| \leq \sum_{j=1}^{l-k} \|C_{j}\|^{2} \|Q^{-1}\| \leq \text{ constant } \sum_{j=1}^{\infty} j^{2\psi} < \infty.
\]

For $k \le l \in \Z_{+}$ let 
\[
  \rho_{\H_{k-1}}(\H_{0,k, \D}, \H_{(l-r)^{+},\rwrev{l-1}, B}) = \sup \big| \Corr_{\H_{{k-1}}}[Z, Y] \big|
\]
denote the maximal conditional correlation coefficient between $\H_{0,k, \D} $ and $\H_{(l-r)^{+},\rwrev{l-1},\D}$ given $\H_{k-1}$, where the supremum is over all random variables $Z$  and $Y$ that are respectively $\H_{0,k, \D}$ and $\H_{(l-r)^{+},\rwrev{l-1}, \D}$ measurable, and where $\EE_{\H_{k-1}}[Z^2]<\infty$ and $\EE_{\H_{k-1}}[Y^2]<\infty$ almost surely.
It follows from~\eqref{eq:condcovboundautoreg} that as $m\to\infty$,
\[
  \sup_{\substack{l,k',k\in\Z_{+}\\ (l-r)^{+} \leq k' \leq l\\ l-k \geq m}} \sup_{\l\in \D,\j\in \D}  \left| \Corr_{\H_{k-1}} \left[ X_{\rwrev{k}}(\l), X_{k'}(\j) \right] \right|\leq \text{ constant } m^{\psi},
\]
which means that
\[
  \sup_{\substack{l,k\in\Z_{+} \\ l-k\geq m}}   \rho_{\H_{k-1}} (\H_{0,k,\D}, \H_{(l-r)^{+},\rwrev{l-1}, \D}) \leq \text{ constant } m^{\psi}.
\]

\subsection{A pseudo-likelihood score function}
\label{sec:pseudologlikARexample}
Given $\H_{0}$, the joint distribution of $X_{1},\ldots,X_{K}$ is given by
\[
  [X_{1}, \ldots, X_{K}]_{\H_{0}} = \prod_{k=1}^{K} [X_{k}]_{\H_{k-1}},
\]
where $[\,\cdot\,]_{\H_{k-1}}$ denotes conditional distribution given $\H_{k-1}$. However, given $\H_0$, a pseudo-likelihood function for $[X_{k}]_{\H_{k-1}}$ is given by
\[
  \prod_{\l\in \D} [X_{k}(\l)]_{\sigma(\H_{k-1} \cup \H_{k, k, \D\backslash\{\l\}})}.
\]
Thus, a log pseudo-likelihood for the observations $X_{1},\ldots, X_{K}$ is given by
\begin{align*}
  \ell(\beta,\gamma) = -\frac{1}{2}\sum_{k=1}^{K} \sum_{\l\in \D} a(\l) \left( X_{k}(\l) - \beta \xi_{k-1}^{\mathrm{temp}}(\l) - \gamma \xi_{k}^{\mathrm{spat}}(\l) \right)^{2}.
\end{align*}
The corresponding score function of $\ell(\beta, \gamma)$ is given by
\[
  T_{K} = \frac{\dd}{\dd (\beta,\gamma)^\T} \ell(\beta,\gamma) = \sum_{k=1}^{K} \sum_{\l\in \D} E_{k}(\l),
\]
where
\[
  E_{k}(\l) = \varepsilon_{k}(\l) \left(\xi_{k-1}^{\mathrm{temp}}(\l), \xi_{k}^{\mathrm{spat}}(\l)\right)^{\T}
\]
and $\varepsilon_{k}=A( X_{k} - \beta \xi_{k-1}^{\mathrm{temp}} - \gamma \xi_{k}^{\mathrm{spat}} )$.

It follows from~M-(6) 
that given $\H_{k-1}$ and with $\beta$ and $\gamma$ the `true' parameter values,
\[
  X_{k} - \beta \xi_{k-1}^{\mathrm{temp}} \sim Q^{-1/2} Z_{k}
\]
and
\[
  \xi_{k}^{\mathrm{spat}} = B_{0} \left( X_{k} - \beta \xi_{k-1}^{\mathrm{temp}} \right) \sim B_{0} Q^{-1/2} Z_{k},
\]
which means that
\begin{align*}
    \varepsilon_{k}=A \left( X_{k} - \beta
  \xi_{k-1}^{\mathrm{temp}} -  \gamma \xi_{k}^{\mathrm{spat}}
  \right) & \sim  A\left( I_{|\D|} - \gamma B_{0} \right) Q^{-1/2} Z_{k} \\
    &= Q Q^{-1/2} Z_{k} = Q^{1/2} Z_{k}.
\end{align*}

\subsection{Covariance matrix of the score function}
\label{sec:covmatrixARexample}

Turning to $\EE_{\H_0}\Var_{H_{k-1}}\sum_{\l \in \D} \Ekn{k}(\l)$,
we have
\[
  \sum_{\l \in \D} \Ekn{k} (\l)= \left ( Z_{k}^\T Q^{1/2} \xi_{k-1}^{\mathrm{temp}}, Z_{k}^\T K Z_{k} \right)^\T,
\]
where $K=[K(\i,\j)]_{\i \j} = Q^{1/2} B_{0}Q^{-1/2}$. Obviously, $\EE_{\H_{k-1}}\left[ Z_{k}^\T Q^{1/2} \xi_{k-1}^{\mathrm{temp}} \right] = 0$ and
\[
  \Var_{\H_{k-1}}\left[ Z_{k}^\T Q^{1/2}  \xi_{k-1}^{\mathrm{temp}} \right] = (\xi_{k-1}^{\mathrm{temp}})^\T
  Q\xi_{k-1}^{\mathrm{temp}}.
\]
Since $Z_{k}^{\T} K Z_{k}$ is a quadratic form in the standard Gaussian random vector $Z_{k}$, we have (see \cite{rencher2008linear}, Theorems 5.2a, 5.2c and 5.2d)
\begin{align*}
  \EE_{\H_{k-1}}\left[ Z_{k}^{\T} K Z_{k}\right] &= \mathrm{trace}(K) = \mathrm{trace}(Q^{1/2} B_{0}Q^{-1/2}) \\
  &=  \mathrm{trace}( B_{0} ) = \sum_{\l\in\D}b_0(\l,\l)=0,\\
 \Var_{\H_{k-1}} Z_{k}^\T K Z_{k} &= 2 \mathrm{trace}(K^{2}) = 2 \mathrm{trace}(Q^{1/2} B_{0} B_{0} Q^{-1/2}) = 2 \mathrm{trace}( B_{0}^{2} )
\end{align*}
and
\[
  \Cov_{\H_{k-1}} \left [ Z_{k}^\T Q^{1/2}
  \xi_{k-1}^{\mathrm{temp}}, Z_{k}^\T K Z_{k} \right ]=2\left(\xi_{k-1}^{\mathrm{temp}} \right)^{\T} Q^{1/2} K \EE_{\H_{k-1}}\left[ Z_{k} \right] = 0.
\]
Thus $T_{K}$ provides an unbiased estimating function and
\[
  \Sigma_{K} = \Var_{\H_{0}}[T_{K}] = \sum_{k=1}^K \EE_{\H_0}\Var_{H_{k-1}}\sum_{\l \in \D} \Ekn{k}(\l) =\mathrm{diag}(\lambda^{(1)}_{K}, \lambda^{(2)}_{K}),
\]
where
\[ \lambda^{(1)}_{K} = \sum_{k=1}^K \EE_{\H_0}\left[ (\xi_{k-1}^{\mathrm{temp}})^\T
  Q\xi_{k-1}^{\mathrm{temp}}\right] \]
and $\lambda_{K}^{(2)} = 2 K  \mathrm{trace}( B_{0}^{2})$.

On the other hand, given $\H_{0}$, from~\eqref{eq:recursive Xkn},
\[
  \beta \xi_{k-1}^{\mathrm{temp}} = X_{k} - Q^{-1/2} Z_{k} = C_{k} X_{0} + \sum_{j=1}^{k-1} C_{j} Q^{-1/2} Z_{k-j}
\]
and hence $\EE_{\H_{0}}\left[\beta \xi_{k-1}^{\mathrm{temp}} \right] = C_{k} X_{0}$ and $\Var_{\H_{0}}\left[ \beta \xi_{k-1}^{\mathrm{temp}} \right] = \sum_{j=1}^{k-1} C_{j} Q^{-1} C_{j}^{\T}$. Thus
\begin{align*}
    \beta^{2} \EE_{\H_0} & \left[ ( \xi_{k-1}^{\mathrm{temp}})^\T   Q\xi_{k-1}^{\mathrm{temp}}\right] \\
    &= \mathrm{trace}\left( Q \Var_{\H_{0}}\left[ \beta \xi_{k-1}^{\mathrm{temp}} \right] \right)
    +  \left( \EE_{\H_{0}}\left[\beta \xi_{k-1}^{\mathrm{temp}} \right] \right)^{\T} Q \EE_{\H_{0}}\left[\beta \xi_{k-1}^{\mathrm{temp}} \right]\\
    &= \sum_{j=1}^{k-1} \mathrm{trace}\left(Q C_{j} Q^{-1} C_{j}^{\T} \right)
    + X_{0}^{\T} C_{k}^{\T} Q C_{k} X_{0}
\end{align*}
and because the precision matrix $Q$ is positive definite,
\begin{align*}
   \beta^{2} \EE_{\H_0} \left[ ( \xi_{k-1}^{\mathrm{temp}})^\T   Q\xi_{k-1}^{\mathrm{temp}}\right] &\geq \sum_{j=1}^{k-1} \mathrm{trace}\left(Q C_{j} Q^{-1} C_{j}^{\T} \right) \\
   &\geq \sum_{j=1}^{k-1} \mathrm{trace}\left(C_{j} C_{j}^{\T} \right) \geq \mathrm{trace}\left(C_{1} C_{1}^{\T} \right).
\end{align*}
But $C_{1} = H^{\T} F H = \beta B_{1}$, so $\mathrm{trace}( C_{1} C_{1}^{\T}) = \beta^{2} \mathrm{trace}\left( B_{1}^{2} \right)$ and $\lambda_{K}^{(1)}\geq K \mathrm{trace} \left( B_{1}^{2} \right)$.

\subsection{Moment bounds}\label{sec:momentbounds}
For any $k\in\N$ and $\nu>0$,
\begin{align*}
  \EE_{\H_{k-1}}\left[ \left| \varepsilon_{k}(\l) \xi_{k-1}^{\mathrm{temp}} (\l) \right|^{\nu} \right] &\leq \left\|\xi_{k-1}^{\mathrm{temp}}  \right\|^{\nu}  \EE_{\H_{k-1}}\left[ \left\| \varepsilon_{k} \right\|^{\nu} \right], \\
  \EE_{\H_{k-1}}\left[ \left| \varepsilon_{k}(\l) \xi_{k}^{\mathrm{spat}} (\l) \right|^{\nu} \right] &\leq \left( \EE_{\H_{k-1}}\left[ \left\| \varepsilon_{k} \right\|^{2\nu} \right] \EE_{\H_{k-1}} \left[ \left\|\xi_{k}^{\mathrm{spat}}  \right\|^{2\nu} \right]  \right)^{\frac{1}{2}}.
\end{align*}
Further,
\begin{align*}
  \EE_{\H_{k-1}}\left[ \left\| \varepsilon_{k} \right\|^{\nu} \right] &= \EE_{\H_{k-1}}\left[ \left\| Q^{1/2}Z_{k} \right\|^{\nu} \right] = \EE_{\H_{k-1}}\left[  \left( Z_{k}^{\T} Q Z_{k} \right)^{ \frac{\nu}{2} }\right], \\
  \EE_{\H_{k-1}}\left[ \left\| \varepsilon_{k} \right\|^{2\nu} \right] &= \EE_{\H_{k-1}}\left[ \left\| Q^{1/2}Z_{k} \right\|^{2\nu} \right] = \EE_{\H_{k-1}}\left[  \left( Z_{k}^{\T} Q Z_{k} \right)^{\nu }\right], \\
  \EE_{\H_{k-1}}\left[ \left\| \xi_{k}^{\mathrm{spat}} \right\|^{2\nu} \right] &= \EE_{\H_{k-1}}\left[ \left( Z_{k}^{\T} Q^{-1/2} B_{0}^{\T} B_{0} Q^{-1/2} Z_{k}  \right)^{\nu} \right]
\end{align*}
are moments of quadratic forms in a standard normal random vector
which are  uniformly bounded in $k$.

Hence,
\begin{align*}
  \EE_{\H_{0}}\left[ \sup_{\l\in \D}\EE_{\H_{k-1}}\left[ \left| \varepsilon_{k}(\l) \xi_{k-1}^{\mathrm{temp}} (\l) \right|^{\nu} \right] \right] &\leq \text{ constant }  \EE_{\H_{0}}\left[ \left\|\xi_{k-1}^{\mathrm{temp}}  \right\|^{\nu}  \right], \\
  \EE_{\H_{0}}\left[ \sup_{\l\in \D}\EE_{\H_{k-1}}\left[ \left| \varepsilon_{k}(\l) \xi_{k}^{\mathrm{spat}} (\l) \right|^{\nu} \right] \right] &\leq \text{ constant.}
\end{align*}
For $\nu >1$, using the representation~\eqref{eq:recursive Xkn} and the
Minkowski inequality,
\begin{align*}
   \EE_{\H_{0}}\left[ \left\| \xi_{k-1}^{\mathrm{temp}} \right\|^{\nu} \right] &= \EE_{\H_{0}}\left[ \left\| C_{k} X_{0} + \sum_{j=1}^{k-1} C_{j} Q^{-1/2} Z_{k-j}  \right\|^{\nu} \right] \\
  &\leq \EE_{\H_{0}}\left[ \left\| C_{k} X_{0} \right\|^{\nu} \right] + \EE_{\H_{0}}\left[ \left\|\sum_{j=1}^{k-1} C_{j} Q^{-1/2} Z_{k-j}  \right\|^{\nu} \right] \\
    & \leq \|C_{k}\|^{\nu} \|X_{0}\|^{\nu} + \sum_{j=1}^{k-1} \|C_{j}\|^{\nu} \EE_{\H_{0}} \left[ \|Q^{-1/2} Z_{k-j}\|^{\nu} \right] \\
    &\leq \text{ constant}_1 k^{\psi\nu} + \text{ constant}_2 \sum_{j=1}^{\infty} j^{\psi\nu} < \infty
\end{align*}
because
\[
  \EE_{\H_{0}} \left[ \|Q^{-1/2} Z_{k-j}\|^{\nu} \right] = \EE_{\H_{0}}\left[ \left( Z_{k-j}^{\T} Q^{-1} Z_{k-j} \right)^{\frac{\nu}{2}} \right]
\]
is a bounded constant that does not depend on $j$ and $k$ and
$C_{k}=O(k^{\psi})$, as $k\to \infty$. Thus, with
\[
  E_{k}(\l) = \varepsilon_{k}(\l) \left(\xi_{k-1}^{\mathrm{temp}}(\l), \xi_{k}^{\mathrm{spat}}(\l)\right)^{\T},
\]
condition M-($\mathcal{A}$4-2) 
holds.



\bibliographystyle{apalike} 
\bibliography{biblio_biscio}